\documentclass[preprint]{elsarticle}
\makeatletter
\def\ps@pprintTitle{%
 \let\@oddhead\@empty
 \let\@evenhead\@empty
 \def\@oddfoot{\centerline{\thepage}}%
 \let\@evenfoot\@oddfoot}
\makeatother
\usepackage[english]{babel}
\usepackage{inputenc} 
\usepackage[T1]{fontenc}
\usepackage{fancyhdr}
\usepackage{graphicx}

\usepackage{amsmath, amsfonts, amssymb, amsthm}

\title{An accurate integral equation method for \\ simulating multi-phase Stokes flow}
\author{Rikard Ojala and Anna-Karin Tornberg}

\address{KTH Mathematics, Linn\'{e} FLOW Centre/Swedish e-Science Research Centre, \\ 100 44 Stockholm, Sweden}

\begin{document}

\begin{abstract}
We introduce a numerical method based on an integral equation formulation for simulating drops in viscous fluids in the plane. It builds upon the method introduced by Kropinski in 2001 \cite{Kropinski2001}, but improves on it by adding an interpolatory quadrature approach for handling near-singular integrals.  Such integrals typically arise when drop boundaries come close to one another, and are difficult to compute accurately using standard quadrature rules. Adapting the interpolatory quadrature method introduced by Helsing and Ojala in 2008 \cite{Helsing2008} to the current application, very general drop configurations can be handled while still maintaining stability and high accuracy. The performance of the new method is demonstrated by some challenging numerical examples.
\end{abstract}

\maketitle
\section{Introduction} 

There is a rapidly growing research area within fluid mechanics referred to as ''micro-fluidics'', which is partly driven by the strong trend to miniaturize equipment for chemical analysis and synthesis.  One emerging technology is ''droplet micro-fluidics'' in which an aqueous sample is emulsified into droplets and dispersed in a continuous oil phase \cite{Joensson2012,Seemann2012}. At these small scales, viscous forces dominate and inertial effects are negligible, and the Reynolds numbers are very small. With pico liter sized droplets, the surface to volume ratio is large, and the interface dynamics is becoming increasingly important.

Over the years, numerous numerical methods have been developed for the
simulation of immiscible multiphase flow.  In a sharp interface
mathematical description of two immiscible fluids on the continuum
level, the fluid-fluid interface is assumed to be infinitesimally
thin. Surface forces are singularly supported on the interfaces
separating the two fluids, leading to a discontinuity in pressure and
velocity gradients across the fluid interfaces. Due to these
difficulties, it still remains a challenge to perform highly accurate
simulations of immiscible multiphase flow. 

The dominant class of methods in the literature is what is often called interface tracking or interface
capturing methods. In these methods, the Navier-Stokes or
Stokes equations are solved on a computational grid or mesh, that is
not required to conform to the fluid-fluid interfaces, and the interfaces are
represented separately. Surveys of the main families of methods can be
found in (Level-set methods)  Sethian and Smereka \cite{Sethian2003},
(Front-tracking) Tryggvason et al.\cite{Tryggvason2001} and (Volume
of Fluid methods) Scardovelli and Zaleski \cite{Scardovelli2003}. 

%In all these methods the effect of the surface tension forces (can be
%computed using the interface representation) must be included in the
%fluid solver, and the velocity for the evolution of the
%interface must be computed using the velocity obtained from the fluid solver. 

The most commonly used approach of including the effect of the singular surface tension
force in the fluid solver is to solve the Navier-Stokes equations
with regularized surface tension force, an idea that was introduced for
elastic interfaces by Peskin already in 1977 \cite{Peskin1977}. This
smears the solution at the interface over a thin region, and
velocities will never be better than first order accurate close to the interface, even
if constructions based on specific choices of regularized delta functions
can allow for higher orders of accuracy away from the interface
\cite{Tornberg2004a}. 

In order to increase the accuracy close to the interface, methods that
avoid regularization and instead directly enforce jump conditions at
the interface have been developed, starting with the Immersed Interface
(IIM) method \cite{Leveque1997} and the extended finite element (XFEM)
method \cite{Fries2010}. These methods can achieve second order
accuracy, but issues remain, such as possible
ill-conditioning in XFEM depending on how interfaces intersect the
underlying grid. This remains an area of active research, with new
promising ideas \cite{Wadbro2013}. 

The Stokes equations for multiphase flow can be reformulated as
boundary integral equations containing integrals over the fluid-fluid
interfaces, see \cite{Pozrikidis1992}. The above mentioned issues are
then avoided - there is no underlying volume grid to couple to, jumps
in solutions are naturally taken care of, and viscosity ratios between
fluids enter only in coefficients of the equations.

When the Stokes equations are reformulated as boundary integral
equations, the dimension of the equations and hence the number of
unknowns in the discretized problem is reduced. The drawback is that
these discretizations result in dense linear systems, making them
very costly to solve.  This can be
addressed by the use of acceleration techniques such as either a fast
multipole (FMM) method \cite{Greengard1997, Wang2007} or FFT-based
methods \cite{Saintillan2005,Lindbo2010,Lindbo2011b}.

Discretizations with a so called Nystr\"{o}m method will be of high
accuracy, given that the errors in the numerical integration are kept
small. The integrands can however be weakly singular or singular and special
quadrature is needed to obtain high accuracy. The case of {\em near
  singularity} occurs when the evaluation point is close to the
interface but not on it, such as when it is located on an interface nearby. 
For the singular case, a large gain in accuracy can be obtained by
local modifications to a quadrature rule in the vicinity of the
singularity \cite{Kapur1997,Marin}, modifications that are specific to the
class of singular functions that is considered. One can also introduce
a mapping that removes the principal singularity, see e.g.  \cite{Bruno2001}.
These methods are not applicable for nearly singular integrals. 
For this case, Helsing and Ojala \cite{Helsing2008} have developed
an interpolatory scheme that offers very high accuracy for the integral
kernels of Laplace's equation in 2D. Another new and promising
development is the so-called QBX method, designed for the 2D Helmholtz
equation in \cite{Kloeckner2013}. 

In this paper we present a numerical scheme for simulating Stokes flow capable of 
handling very general drop configurations while maintaining high accuracy. We will use
the integral equation formulation from Kropinski, \cite{Kropinski2001}, and by adapting the 
interpolatory quadrature introduced in \cite{Helsing2008}, we are able 
to allow drops to be arbitrarily close to each other without losing stability or accuracy. This 
is the main novelty of the paper, and opens up for drop simulations that were previously infeasible.

The structure of this paper is as follows : in section \ref{sec:problem} we discuss the specific
problem we are solving and introduce its integral equation formulation. In section \ref{sec:quad}, we
describe the interpolatory quadrature approach we will use, and in section \ref{sec:numm} we discuss the
numerical method in some detail. In section \ref{sec:num} we demonstrate how the solver performs 
for some drop configurations, and lastly, in \ref{sec:outlook} some options for future extensions of the solver are outlined.

\section{Problem statement}
\label{sec:problem}
\begin{figure}[h]
\centering
\resizebox{65mm}{!}{\includegraphics{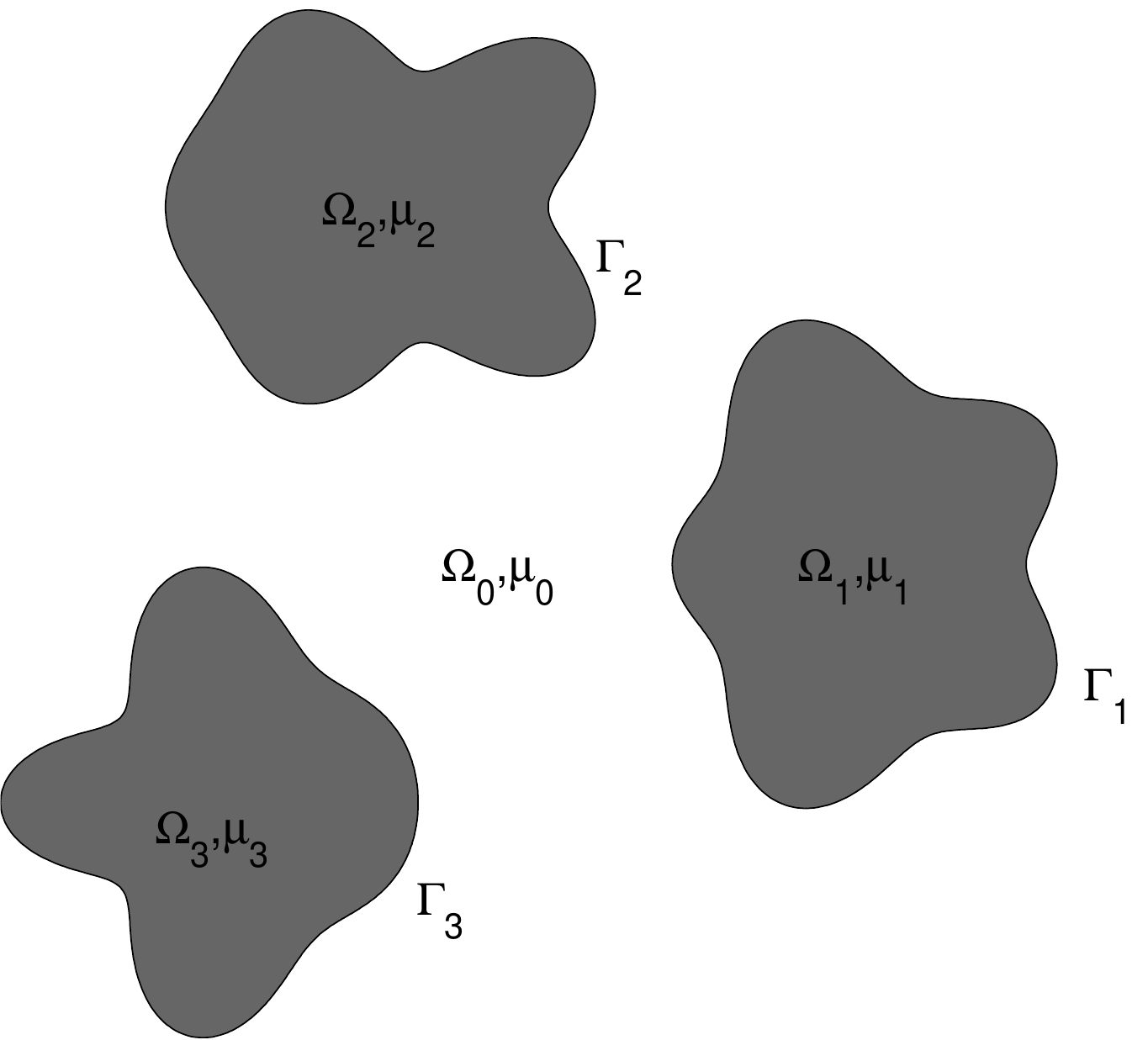}}
\caption{A drop configuration with three drops.}
\label{fig:slanted}
\end{figure}
Consider an infinite expanse of fluid in the plane with viscosity $\mu_0$ containing $n$ drops $\Omega_k$ with boundary $\Gamma_k$ and viscosity $\mu_k = \lambda_k\mu_0$ for $k = 1,2,...,n$. The surrounding fluid resides in $\Omega_0$. All viscosities are positive and we set $\Gamma$ to be the union of all $\Gamma_k$. Assuming zero Reynolds number flow, the equations governing the flow are the Stokes equations 
\begin{eqnarray}
\mu_0 \Delta \mathbf{u}_0 &=& \nabla p_0, \quad \nabla \cdot \mathbf{u}_0 = 0, \quad \mathbf{x} \in \Omega_0, \label{eq:stokes1}\\
\mu_k \Delta \mathbf{u}_k &=& \nabla p_k, \quad \nabla \cdot \mathbf{u}_k= 0, \quad \mathbf{x} \in \Omega_k, \label{eq:stokes2}
\end{eqnarray}
where $\mathbf{u}_0$ and $p_0$ is the velocity and pressure in the surrounding fluid and $\mathbf{u}_k$ and $p_k$ is the velocity and pressure in drop $k$. In this paper, the flow will be driven by surface tension only, and we will assume that the $n$ drops will not be in contact with each other, only with the surrounding fluid. The boundary conditions at the drop boundaries are therefore continuity of the velocity and a jump in normal stress proportional to the curvature. On the boundary of drop $k$ we write this latter condition as
\begin{equation}
-(p_0-p_k)\mathbf{n} + 2(\mu_0\boldsymbol{\varepsilon}_0 - \mu_k\boldsymbol{\varepsilon}_k)\mathbf{n} = -\sigma\kappa\mathbf{n}, \label{eq:bndrycond}
\end{equation}
where $\mathbf{n}$ is the outward unit normal, $\sigma$ is the surface tension coefficient, $\kappa$ is the curvature and $\boldsymbol{\varepsilon}$ is the $2\times 2$ rate of strain tensor with elements 
\begin{equation}
\varepsilon_{ij} = \frac{1}{2}\left(\frac{\partial u_i}{\partial x_j} + \frac{\partial u_j}{\partial x_i}\right). 
\end{equation}
We use the same surface tension coefficient $\sigma$ for all drop interfaces. It is convenient to non-dimensionalize the equations to get rid of a few of the physical parameters. We do this by setting
\begin{eqnarray}
\tilde{\mathbf{u}}_k &=& \frac{\mu_0}{\sigma}\mathbf{u}_k, \\
\tilde{p}_k &=& \frac{L}{\sigma}p_k,\\
\tilde{\mathbf{x}} &=&\frac{1}{L}\mathbf{x},
\end{eqnarray}
where $L$ is the characteristic length which we take to be the length of the computational domain encompassing $\Gamma$. To simplify notation we do away with the tildes and just look at the non-dimensionalized equations. The equations (\ref{eq:stokes1}) and (\ref{eq:stokes2}) now become
\begin{eqnarray}
\Delta \mathbf{u}_0 &=& \nabla p_0, \quad \nabla \cdot \mathbf{u}_0 = 0, \quad \mathbf{x} \in \Omega_0, \label{eq:stokes3}\\
\lambda_k \Delta \mathbf{u}_k &=& \nabla p_k, \quad \nabla \cdot \mathbf{u}_k= 0, \quad \mathbf{x} \in \Omega_k, \quad k=1,2,...,n \label{eq:stokes4}
\end{eqnarray}
and the normal stress boundary condition (\ref{eq:bndrycond}) becomes
\begin{equation}
-(p_0-p_k)\mathbf{n} + 2(\boldsymbol{\varepsilon}_0 - \lambda_k\boldsymbol{\varepsilon}_k)\mathbf{n} = -\kappa\mathbf{n},
\end{equation}
which is enforced on each $\Gamma_k, k=1,2,...,n$. We will assume that $\mu_0 > 0$ and $\lambda_k > 0$, but it should be remarked that the limiting cases $\lambda_k = 0$ and $\mu_0 = 0$, corresponding to inviscid bubbles and an inviscid surrounding medium respectively, can be handled by the methods in this paper with some adjustments. See \cite{Kropinski2001} and \cite{Kropinski2002} for details.

In order to propagate the drop boundaries forward in time we use a quasi-static approach. That is, we compute the fluid velocity $\mathbf{u}(\mathbf{x})$ on $\Gamma$ via the steady Stokes equations (\ref{eq:stokes3}) and (\ref{eq:stokes4}) and propagate $\Gamma$ using an ODE. We thus need a non-dimensionalized time variable which is $\tilde{t} = \sigma/(L\mu_0) t$. The natural Lagrangian approach
\begin{equation}
\frac{\partial \mathbf{x}}{\partial t} = \mathbf{u}(\mathbf{x}), \qquad \mathbf{x} \in \Gamma, \label{eq:ODE}
\end{equation}
is not suitable as it stands, as explained in for example \cite{Hou1994} and \cite {Kropinski2001}, since the points making up the interface tend to cluster in some areas and thin out in others. The clustering causes the ODE to become stiff and the thinning out can make for insufficient resolution of the interface. It has been pointed out elsewhere that only the normal velocity affect the shape of the interface: the tangential component can be chosen as to preserve the distribution of the points of the interface dynamically. That is, we solve an ODE such as (\ref{eq:ODE}), but with a modified $\mathbf{u}(\mathbf{x})$. We will discuss the implementation of this in section \ref{sec:arcl}.

To solve (\ref{eq:stokes3})-(\ref{eq:stokes4}) and compute the velocity on $\Gamma$, we will use the Sherman-Lauricella integral equation formulation for Stokes flow presented in \cite{Kropinski2001}. In short, it can be shown that the Stokes equations can be reduced to the biharmonic equation, which, in turn, can be recast as a problem in the theory of analytic functions. We express the solution using two analytic functions called Goursat functions, and by choosing the form of the Goursat functions in a clever way we both fulfill the continuity of the velocity and get a resulting Fredholm integral equation that is of the second kind with compact integral kernels. Such equations are often very well conditioned, promising high accuracy solutions that can be computed rapidly using iterative linear equation solvers. 

In what follows, we will be working in the complex plane. We reserve $z$ and $\tau$ to be points in the plane, and $\Re e\left\{f(z)\right\}$ and $\Im m\left\{f(z)\right\}$ denotes the real and and imaginary parts of the complex function $f(z)$. Defining $\lambda(z) = \lambda_k$ for $z \in \Gamma_k$, and setting $\beta(z) = \frac{1-\lambda(z)}{1+\lambda(z)}$ and $\gamma(z) = \frac{1}{1+\lambda(z)}$, the resulting integral equation to solve is
\begin{multline}
\omega(z) + \frac{\beta(z)}{\pi}\int_\Gamma \omega(\tau)\Im m \left\{\frac{{\rm d}\tau}{\tau -z}\right\} + \frac{\beta(z)}{\pi}\int_\Gamma \overline{\omega(\tau)} \frac{\Im m \left\{{\rm d} \tau(\bar{\tau}-\bar{z})\right\}}{(\bar{\tau}-\bar{z})^2} + \\ + \beta(z)\int_\Gamma \omega(\tau)|{\rm d}\tau| = - \frac{\gamma(z)}{2}\frac{\partial z}{\partial s}, \qquad z \in \Gamma, \label{eq:inteq}
\end{multline}
for the unknown complex density $\omega(z)$, where we have set the far-field reference pressure to zero. Here, bar means complex conjugation and by $\frac{\partial z}{\partial s}$ we mean differentiation with respect to arclength. The integral operators are both compact on smooth $\Gamma$. In fact, the integral kernels are bounded; limits exist for $\tau = z$. The first two integrals in (\ref{eq:inteq}) are actually the complex variable version of the integral over the two dimensional stresslet, more commonly written in real variables as
\begin{equation}
\int_\Gamma \omega_i(\mathbf{x}) T_{ijk}(\mathbf{x},\mathbf{x_0}) n_k(\mathbf{x}) {\rm d}s(\mathbf{x}), \qquad j = 1,2,
\end{equation}
where the integration is with respect to arclength and where
\begin{equation}
T_{ijk}(\mathbf{x},\mathbf{x_0}) = -4\frac{\bar{x}_i\bar{x}_j\bar{x}_k}{r^4}, \qquad \bar{\mathbf{x}} = \mathbf{x}-\mathbf{x_0}, \qquad r = |\mathbf{x}|.
\end{equation}
 
Note that if $\lambda_k = 1$ for all $k$, that is if the viscosity of all the drops are the same as for the surrounding fluid, then $\beta(z)=0$, and the solution of the integral equation (\ref{eq:inteq}) is trivial. This is a characteristic shared with primitive variable formulations of the same problem \cite{Pozrikidis1992}, and gives opportunities for rapid testing. 

Having computed $\omega(z)$ the velocity on $\Gamma$ is calculated via
\begin{equation}
u_1(z) + {\rm i}u_2(z) = -\frac{1}{\pi}\int_\Gamma\omega(\tau)\Re e\left\{\frac{{\rm d}\tau}{\tau -z}\right\} - \frac{1}{\pi {\rm i}}\int_\Gamma \overline{\omega(\tau)} \frac{\Im m \left\{{\rm d} \tau(\bar{\tau}-\bar{z})\right\}}{(\bar{\tau}-\bar{z})^2}, \qquad z \in \Gamma, \label{eq:vel}
\end{equation}
where the first integral is singular and must be interpreted in a principal value sense. Here, $u_1$ and $u_2$ are the $x$ and $y$ components of the velocity. In fact, this formula can be used to compute the velocity of the flow for $z$ in the entire complex plane. We will do so in section \ref{sec:num} when we show flow fields for some drop configurations. It is not hard to add ambient flows such as extensional or shear flows if needed. It involves adding terms to the right hand sides of (\ref{eq:inteq}) and (\ref{eq:vel}), but we will not do so here. See \cite{Kropinski2001} for details. 

It has been mentioned above that the integral equation (\ref{eq:inteq}) is well-behaved, and even though the formula (\ref{eq:vel}) contains a singular integral operator evaluating it using standard high order quadrature rules is not a problem. On the whole, things look straight forward. However, when two drop boundaries are close or when a drop boundary falls back on itself (see for example the C-shape boundary in section \ref{sec:cshape}), evaluation of the integral operators cease to be a trivial task. We will now discuss how to treat the above equations numerically to overcome this problem.

\section{Specialized quadratures}
\label{sec:quad}
A major problem plaguing integral equation based solvers of interface problems is that for close-lying interfaces the resulting velocities become inaccurate. This is because the kernels integrated over are near-singular in such cases and standard quadratures do not work well. In turn, phenomena such as lubrication are not captured correctly. In order to get accurate results, we need to replace the standard quadrature schemes with specialized ones in areas where two or more interfaces are close to each other. The same problem arises when we wish to evaluate the velocity in one of the fluids, not on the interface but close to it. The approach we will use is an adaptation of the interpolatory quadrature scheme in \cite{Helsing2008} for use with the integral kernels arising from the equations of Stokes flow.

The special quadrature we will employ is a local one and works on point-panel pairs. It will be invoked whenever a target point $z$ is close enough to a quadrature panel $\Gamma_p$ to warrant special quadrature treatment, that is, when standard 16-point Gauss-Legendre quadrature does not give accurate results for integration on that panel. Because of this, we will restrict ourselves to integrals over one quadrature panel $\Gamma_p$ rather than the entire boundary $\Gamma$ in this section.

There are two instances where accurate integration is required: when solving the integral equation (\ref{eq:inteq}) and when computing the velocity via (\ref{eq:vel}). The integral
\[
\int_{\Gamma_p} \overline{\omega(\tau)} \frac{\Im m\left\{\left(\bar{\tau}-\bar{z}\right){\rm d}\tau\right\}}{(\bar{\tau}-\bar{z})^2}
\]
is common to them both. It can be rewritten as
\begin{equation}
\int_{\Gamma_p} \overline{\omega(\tau)} \frac{\Im m\left\{\left(\bar{\tau}-\bar{z}\right){\rm d}\tau\right\}}{(\bar{\tau}-\bar{z})^2} = \frac{1}{2{\rm i}} \overline{\int_{\Gamma_p} \frac{\omega(\tau)\overline{n_\tau}^2{\rm d}\tau}{\tau-z}} - \frac{1}{2{\rm i}}\overline{\int_{\Gamma_p} \frac{\omega(\tau)\left(\bar{\tau}-\bar{z}\right){\rm d}\tau}{(\tau-z)^2}}, \label{eq:decomp}
\end{equation}
where $n_\tau$ is the outward unit normal at $\tau$. The two other integrals we need to evaluate are
\begin{equation}
\int_{\Gamma_p} \omega(\tau)\Im m\left\{\frac{{\rm d}\tau}{\tau-z}\right\} \quad {\rm and} \quad \int_{\Gamma_p} \omega(\tau)\Re e\left\{\frac{{\rm d}\tau}{\tau-z}\right\}, \label{eq:imre}
\end{equation}
and splitting the complex density $\omega(\tau)$ into its real and imaginary components for these two integrals we see that all the integrals we will deal with are on the forms
\begin{eqnarray}
I_1 &=& \int_{\Gamma_p} \frac{f(\tau){\rm d}\tau}{\tau - z}, \label{eq:i1}\\
I_2 &=&  \int_{\Gamma_p} \frac{f(\tau){\rm d}\tau}{(\tau - z)^2},\label{eq:i2}
\end{eqnarray}
for some function $f(\tau)$. To compute $I_1$ and $I_2$, we approximate $f(\tau)$ by
\begin{equation}
f(\tau) \approx \sum_{k=0}^{15} c_k\tau^k \label{eq:approx}
\end{equation}
over the panel $\Gamma_p$. Determining the coefficients $c_k$ involves solving a Vandermonde system, but through complex scaling and translation we may assume that the edges of $\Gamma_p$ is at -1 and 1 in the complex plane, keeping the conditioning under control, see Appendix A in \cite{Helsing2008}. Note that this requires a multiplication of a scaling factor for (\ref{eq:i2}), because of the lack of scale invariance of the integrand. It should be noted that for the polynomial approximation of $f(\tau)$ to be accurate, $f(\tau)$ needs to be in some sense well-behaved. In practice, this is seen to be the case, at least as long as the density $\omega(\tau)$ and the boundary is reasonably well resolved. 

Inserting (\ref{eq:approx}) into (\ref{eq:i1}) and (\ref{eq:i2}) gives
\begin{eqnarray}
I_1 &=& \int_{\Gamma_p} \frac{f(\tau){\rm d}\tau}{\tau - z} \approx \sum_{k=0}^{15} c_k \int_{-1}^1 \frac{\tau^k {\rm d}\tau}{\tau-z_0} = \sum_{k=0}^{15} c_k p_k, \label{eq:pk} \\
I_2 &=& \int_{\Gamma_p} \frac{f(\tau){\rm d}\tau}{(\tau - z)^2} \approx \alpha \sum_{k=0}^{15} c_k \int_{-1}^1 \frac{\tau^k {\rm d}\tau}{(\tau-z_0)^2} = \alpha \sum_{k=0}^{15} c_k q_k,
\end{eqnarray}
where $\alpha = (z_2-z_1)/2$ is the scaling factor and $z_1$ and $z_2$ are the endpoints of the untransformed $\Gamma_p$, and $z_0$ is the transformed $z$, rotated and scaled along with $\Gamma_p$. The complex numbers $p_k$ and $q_k$ can be computed rapidly using recursions. For $p_0$ we have
\[
p_0 = \int_{-1}^1 \frac{{\rm d}\tau}{\tau - z_0} = \log(1-z_0)-\log(-1-z_0),
\]
but care must be taken if $z_0$ is between the transformed $\Gamma_p$ and the real axis, in which case the residue $2\pi i$ must be added or subtracted from $p_0$, depending on the orientation of $\Gamma_p$. The rest of the $p_k$'s are then computed via
\[
p_{k} = z_0p_{k-1} + \frac{1-(-1)^k}{k}, \qquad k = 1,2,...,15.
\]
For $q_k$ the situation is similar. We get
\[
q_0 = \int_{-1}^1 \frac{{\rm d}\tau}{(\tau - z_0)^2} = -\frac{1}{1+z_0}-\frac{1}{1-z_0},
\]
and
\[
q_{k} = z_0q_{k-1} + p_{k-1}, \qquad k = 1,2,...,15.
\]
To be more precise, for the first integral in the right hand side of (\ref{eq:decomp}) we use $I_1$ with $f(\tau) = \omega(\tau)\overline{n_\tau}^ 2$ and conjugate the result. For the second we use $I_2$ with $f(\tau) = \omega(\tau)(\bar{\tau}-\bar{z})$ and again conjugate the result. For the integrals in (\ref{eq:imre}) we use $I_1$ with $f(\tau) = \omega(\tau)$ but take the imaginary and real part of $p_k$ respectively when computing the sum in (\ref{eq:pk}).

This special quadrature is invoked for any point-panel pair when two conditions are fulfilled: first, if the target point is within one panel length of the midpoint of the panel and second, if the analytical and numerical calculation of $p_0$ differs by more than some set tolerance. We use $10^{-13}$.
Checking all point-panel pairs would be expensive, so a simple grid structure is set up reducing the number of point-panel checks to the closest few.

It should be noted that the special quadrature approach presented here is in no way restricted to the particular complex variable integral equation formulation used in this paper. It could be used for any formulation provided that the integral operators acting on monomials can be evaluated analytically. This is true for operators on the forms
\[
\int_\Gamma f(\tau) \log(|\tau-z|) {\rm d}\tau \quad {\rm or} \quad \int_\Gamma \frac{f(\tau){\rm d}\tau}{(\tau-z)^m},
\]
where $m = 1,2,...$, which includes primitive variable formulations.

\section{Numerical method}
\label{sec:numm}
\begin{figure}[h]
\centering
\resizebox{105mm}{!}{\includegraphics{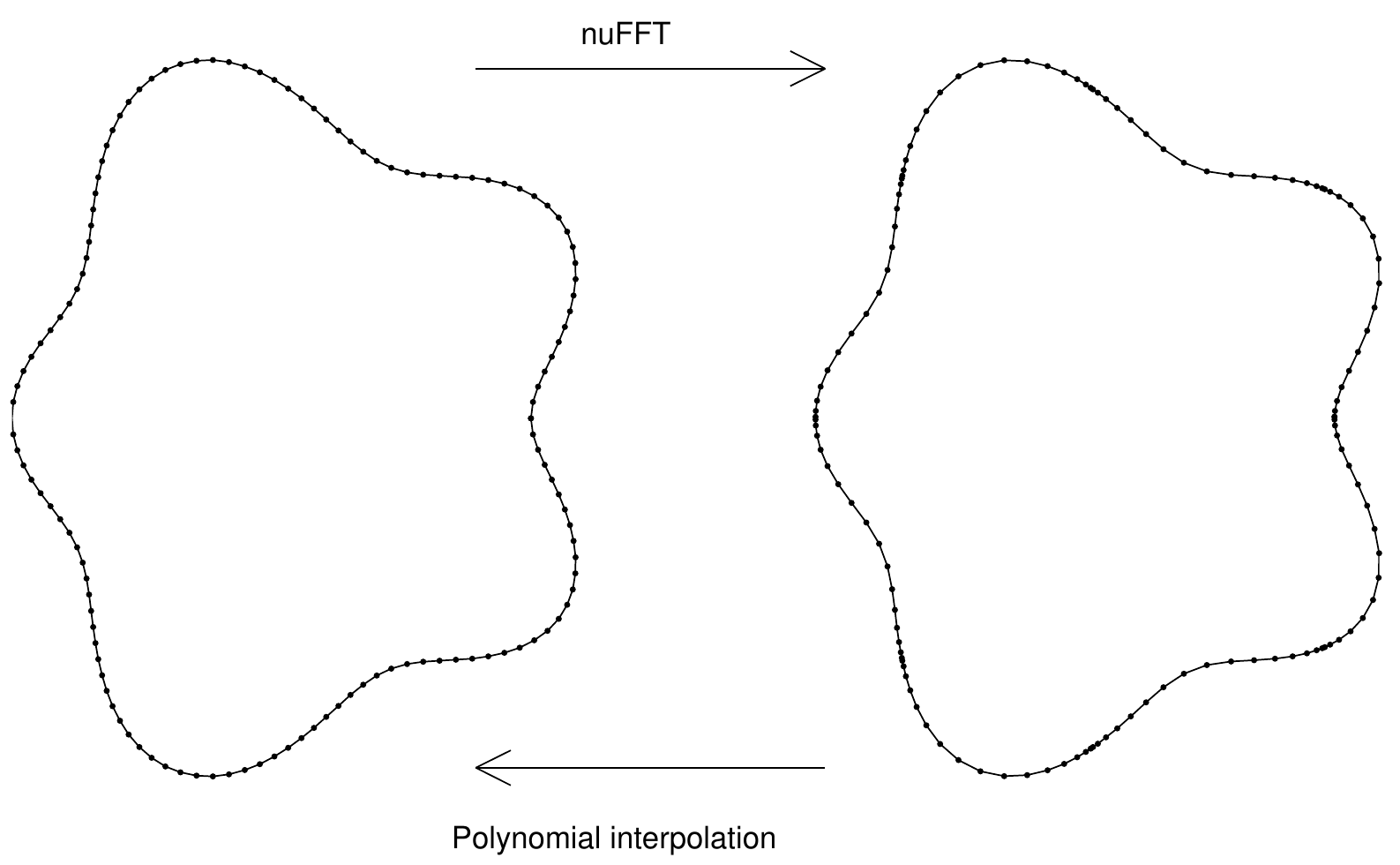}}
\caption{The two grids. Equidistant to the left, and Gauss-Legendre to the right. The non-uniform FFT and polynomial interpolation is used to switch between the grids.}
\label{fig:dropdisc}
\end{figure}
In this section we will discuss the numerical solver, and the components it comprises. 
Suppose we are given $n$ functions $z_k(s), s \in [0,2\pi]$, each describing the initial boundary $\Gamma_k$ of one drop. The drops are not allowed to touch or intersect. We distribute $N_k$ points on each drop boundary $\Gamma_k$ so that they are equispaced in arclength. The equal arclength distribution of points is accomplished by solving a non-linear system of $N_k$ equations using Newton iteration. We end up with $N = \sum_{k=1}^n N_k$ discretization points or marker points describing the boundaries of our drops.

We will be working on two grids, namely the equidistant (trapezoidal rule) grid described above, and a 16-point composite Gauss-Legendre grid, see figure \ref{fig:dropdisc}. The reason for this is twofold:
\begin{enumerate}
\item The scheme preventing the unnecessary clustering of grid points works best on equidistant grids, so the modification of the tangential velocity as well as the time-stepping will be performed on that grid, see section \ref{sec:arcl}.
\item The special quadrature tool we use requires a composite Gaussian grid. The integral equation (\ref{eq:inteq}) will be solved and the velocity computed via (\ref{eq:vel}) on that grid.
\end{enumerate}
Interpolation from the equidistant grid to the Gauss-Legendre grid is done using the non-uniform FFT(nuFFT) with fast Gaussian gridding \cite{Greengard2004} giving high order accuracy in $O(N\log N)$ operations. In the other direction we use standard $O(N)$ 16th degree polynomial interpolation. We keep each $N_k$ a multiple of 16 to facilitate simple interpolation from the Gauss-Legendre grid to the equidistant grid.

When the velocity at the drop boundaries are to be computed we proceed as follows :
\begin{enumerate}
\item Double the number of discretization points using the FFT. Interpolate to a composite 16-point Gauss-Legendre grid using the nuFFT.
\item Solve the discretized version of the integral equation (\ref{eq:inteq}).
\item Compute the velocity at the boundaries from the discretized version of (\ref{eq:vel}).
\item Interpolate back to the equispaced grid using standard 16-degree interpolation and modify the tangential velocity to preserve arclength spacing. Halve the number of points, again using the FFT.
\end{enumerate}
Once the modified velocity on the boundaries have been computed, it is fed into a time-stepper to propagate the boundaries forward in time.

These steps will now be described in some detail in sections \ref{sec:GL-ip}-\ref{sec:arcl}.
To simplify notation we restrict ourselves to one drop in what follows, extending to more than one is straight-forward. 
\subsection{Interpolation to the Gauss-Legendre grid}
\label{sec:GL-ip} 
We start out by temporarily doubling the number of discretization points by taking the FFT of the boundary point vector and padding the Fourier spectrum. This is important for keeping the scheme stable, as the high frequency components tend to grow spuriously \cite{Kropinski2001}. At the end of step 4 the number of points is halved, removing the high-frequency components. We interpolate from the equispaced grid to the Gauss-Legendre grid using the nuFFT, in the process computing derivatives with respect to parameter. We denote by $\tau_i$, $\tau'_i$ and $\tau''_i$ the $2N$ points and derivatives on the boundary. In the context of the Gauss-Legendre grid, we call each group of 16 points a quadrature panel.

\subsection{Solving the integral equation}
\label{sec:solv}
We use the Nystr\"om method to solve the integral equation (\ref{eq:inteq}). The discretized equation reads
\begin{equation}
\omega_i + \frac{\beta_i}{\pi}\sum_{j=1}^{2N} \omega_j M^{(1)}_{ij} + \frac{\beta_i}{\pi}\sum_{j=1}^{2N} \bar{\omega}_j M^{(2)}_{ij} +  \beta_i\sum_{j=1}^{2N}\omega_j|{\rm w}_j\tau_j'| = -\frac{\gamma_i}{2}\frac{\tau_i'}{|\tau_i'|}, \label{eq:disceq}
\end{equation}
for $i = 1,2,...,2N$. Here, $\omega_i$, $\beta_i$ and $\gamma_i$ are the values of the functions $\omega(z)$, $\beta(z)$ and $\gamma(z)$ at $\tau_i$  while ${\rm w}_i$ are the Gauss-Legendre weights. The elements of the matrices $\mathbf{M}^{(1)}$ and $\mathbf{M}^{(2)}$ are
\begin{eqnarray}
M^{(1)}_{ij} &=& \Im m\left\{\frac{{\rm w}_j'\tau_j'}{\tau_j-\tau_i}\right\},\\
M^{(2)}_{ij} &=& \frac{\Im m \left\{{\rm w}_j\tau'_j(\bar{\tau}_j-\bar{\tau}_i)\right\}}{(\bar{\tau}_j-\bar{\tau}_i)^2}.
\end{eqnarray}
For the diagonal elements, where $i=j$, limits are available. These are 
\begin{eqnarray}
M^{(1)}_{ii} &=& \Im m \left\{ \frac{{\rm w}_i\tau''_i}{2\tau'_i}\right\},\\
M^{(2)}_{ii} &=& \frac{\Im m\left\{{\rm w}_i\tau''_i\bar{\tau}'_i\right\}}{2(\bar{\tau}'_i)^2}.
\end{eqnarray}
The system of equations is solved using GMRES \cite{Saad1986} and the fast multipole method \cite{Greengard1987} is used to compute the action of $\mathbf{M}^{(1)}$ and $\mathbf{M}^{(2)}$ rapidly. Owing to the fact that the integral equation is of Fredholm's second kind, the spectral properties of the system (\ref{eq:disceq}) are such that the number of GMRES iterations is bounded for increasing $N$. The number of iterations is determined by the geometry and the viscosity ratios, i.e. by the underlying problem and not the discretization. For complicated drop configurations, or when the viscosity ratios differ greatly, the number of GMRES iterations needed to solve (\ref{eq:disceq}) can be quite high, at least in the initial stages of a simulation. In these cases, direct solvers \cite{Martinsson2005} could be used to speed up the process.

If the drop boundary falls back on itself, or if two or more drops are close to each other, the composite 16-point Gauss-Legendre quadrature fails to be accurate. The reason for this is that the kernels become near-singular. We compute special quadrature corrections to the matrices when this happens, see section \ref{sec:quad} for details. 
\subsection{Computing the velocity}
\label{sec:vel}
Once the discrete value of the density $\omega_i$ at the discrete points $\tau_i$ are available we compute the velocity at $\tau_i$ by a discretized version of (\ref{eq:vel}). Evaluating the velocity is a bit more involved than solving the integral equation above since here we need to deal with a singular integral operator. For convenience we restate the formula :
\begin{equation}
u_1(z) + {\rm i}u_2(z) = -\frac{1}{\pi}\int_\Gamma\omega(\tau)\Re e\left\{\frac{{\rm d}\tau}{\tau -z}\right\} - \frac{1}{\pi{\rm i}}\int_\Gamma \overline{\omega(\tau)} \frac{\Im m \left\{{\rm d} \tau(\bar{\tau}-\bar{z})\right\}}{(\bar{\tau}-\bar{z})^2}, \qquad z \in \Gamma. 
\end{equation}
We note that the second integral is the same as the second one above, and its discretization is given by $\mathbf{M}^{(2)}$. For the first one we use singularity subtraction, and put together we get
\begin{equation}
(u_1)_i + {\rm i}(u_2)_i = - \frac{{\rm w}_i}{\pi}\omega'_i - \frac{1}{\pi}\underset{j\neq i}{\sum_{j=1}^{2N}} (\omega_j-\omega_i) \Re e \left\{\frac{{\rm w}_j\tau'_j}{\tau_j-\tau_i}\right\} - \frac{1}{\pi{\rm i}} \sum_{j=1}^{2N} M_{ij}^{(2)}\bar{\omega}_j
\end{equation}
for the velocity at $\tau_i$, where $\omega'_i$ is the derivative of $\omega_i$ calculated numerically using 16-point interpolation, and ${\rm w}_i$ are again the Gauss-Legendre weights.
 Again, the fast multipole method is used to evaluate the sums. In cases when drop boundaries are close to each other, special quadrature treatment is required here for the same reasons as when solving the integral equation as discussed above. See section \ref{sec:quad}.

\subsection{Preserving arclength distance}
\label{sec:arcl}
Having computed the velocity on the boundary we now wish to modify its tangential component to dynamically preserve the spacing of the boundary points. The derivation of this modification can be found in \cite{Kropinski2001}. Here, we only state results, rewritten to better fit our notation. In the solver, we begin this step by returning to the equispaced grid, that is we interpolate the velocity to equispaced points using 16th order polynomial interpolation on each quadrature panel.

In the following we will assume that the boundary and velocity is parametrized in $s \in [0,2\pi]$ by $z(s)$ and that the velocity $u(s) = u_1(s) + {\rm i}u_2(s)$. As was stated before, the normal component of the velocity decides the evolution of the drop boundary, so we decompose $u(s)$ into its normal and tangential components :
\begin{equation}
u(s) = (u_n(s) + {\rm i}u_t(s))n(s),
\end{equation}
where $n(s)$ is the complex unit outward normal function. We have that $u_n(s) = \Re e \{u(s)\bar{n}(s)\}$ and if $u_t(s)$ is chosen as
\begin{equation}
u_t(s) = \frac{s}{2\pi}\int_0^{2\pi} \Im m\left\{\frac{z''(q)}{z'(q)}\right\}u_n(q) {\rm d}q - \int_0^s \Im m\left\{\frac{z''(q)}{z'(q)}\right\}u_n(q) {\rm d}q, \label{eq:maintain}
\end{equation}
then the equispaced grid will be preserved. Since we are now working on the equidistant grid, using the trapezoidal rule the integrals are evaluated with spectral accuracy. The derivatives and the antiderivative are computed using FFTs. 

The above formula maintains the relative spacing of the equidistant grid. Although it seems to be analytically possible to preserve a Gauss-Legendre grid dynamically, numerically this is not a good idea. The formula (\ref{eq:maintain}) is altered, requiring the adding of terms on a panel-by-panel basis that are poorly resolved by polynomials which causes severe numerical difficulties. This is the main motivation for switching between equispaced and Gauss-Legendre grids, instead of working on Gauss-Legendre grids only.

\subsection{Time stepping and adaptivity}
We use the Bogacki-Shampine embedded second-third order Runge-Kutta method for time stepping. Experiments show that higher order methods do not pay off in terms of the number of velocity computations versus step size. On the interface of each drop, we have a grid that is initially equidistant in arclength, and the algorithm described in section \ref{sec:arcl} maintains the relative arclength distance. As the circumference of the interfaces changes with time, we adaptively change the number of discretization points to keep the actual arclength distance between points close to the initial one. At regular intervals we therefore check how much the circumference of the drop has changed and decrease or increase the number of discretization points accordingly, but only in steps of 16 points to keep the number a multiple of 16. The underlying assumption is that the arclength distance between discretization points given at the start of the simulation is enough to properly resolve the drop throughout the simulation, which may not be true under all circumstances. However, in the numerical experiments performed in this paper this approach to adaptivity works well, and it typically saves around 30\%-50\% on discretization points towards the end of the simulations.

\section{Numerical examples}
\label{sec:num}
In this section we will present a number of numerical examples showing the capabilities of the solver described above. For the convenience of the readers implementing their own solvers, we will for the first two drop configurations tested report some quantities that may serve as reference values for benchmark tests.  Since the steady-state shape of the drops are circles, we have settled for the position of the centers of these circles, estimated when the maximum deviation from a circle is less than $10^{-3}$ in relative terms. More precisely, introducing $\tilde{\tau}_i = \tau_i - \tau_g$, where $\tau_g$ is the center of gravity of a drop, we say that a drop is circular enough when 
\[
r_{{\rm dev}} = \left\| 1-\frac{|\tilde{\tau}_i|}{|\tilde{\tau}|_{\rm mean}} \right\|_\infty < 10^{-3}.
\]
Here, $|\tilde{\tau}|_{\rm mean}$ is the mean of all $|\tilde{\tau}_i|$. For more than one drop, we require that $r_{{\rm dev}}$ for all drops fall below $10^{-3}$.
We also report the approximate simulation time, $t_{\rm steady}$, when this state is reached as well as the area error, $A_{\rm err}$.

For all experiments, the Runge-Kutta and GMRES tolerances are set to $10^{-8}$ and $10^{-10}$, respectively. The computer on which the experiments are run is equipped with a quad core 3.4 GHz Intel Core i7 processor and 8 Gb of RAM. The fast multipole routines as well as the special quadrature computation takes advantage of multiple processor cores, these are also the most time-critical components in the solver. These components are written in C++, the rest of the solver is written in Matlab.
 
\subsection{The flower}
\begin{figure}[h]
\centering
\resizebox{85mm}{!}{\includegraphics{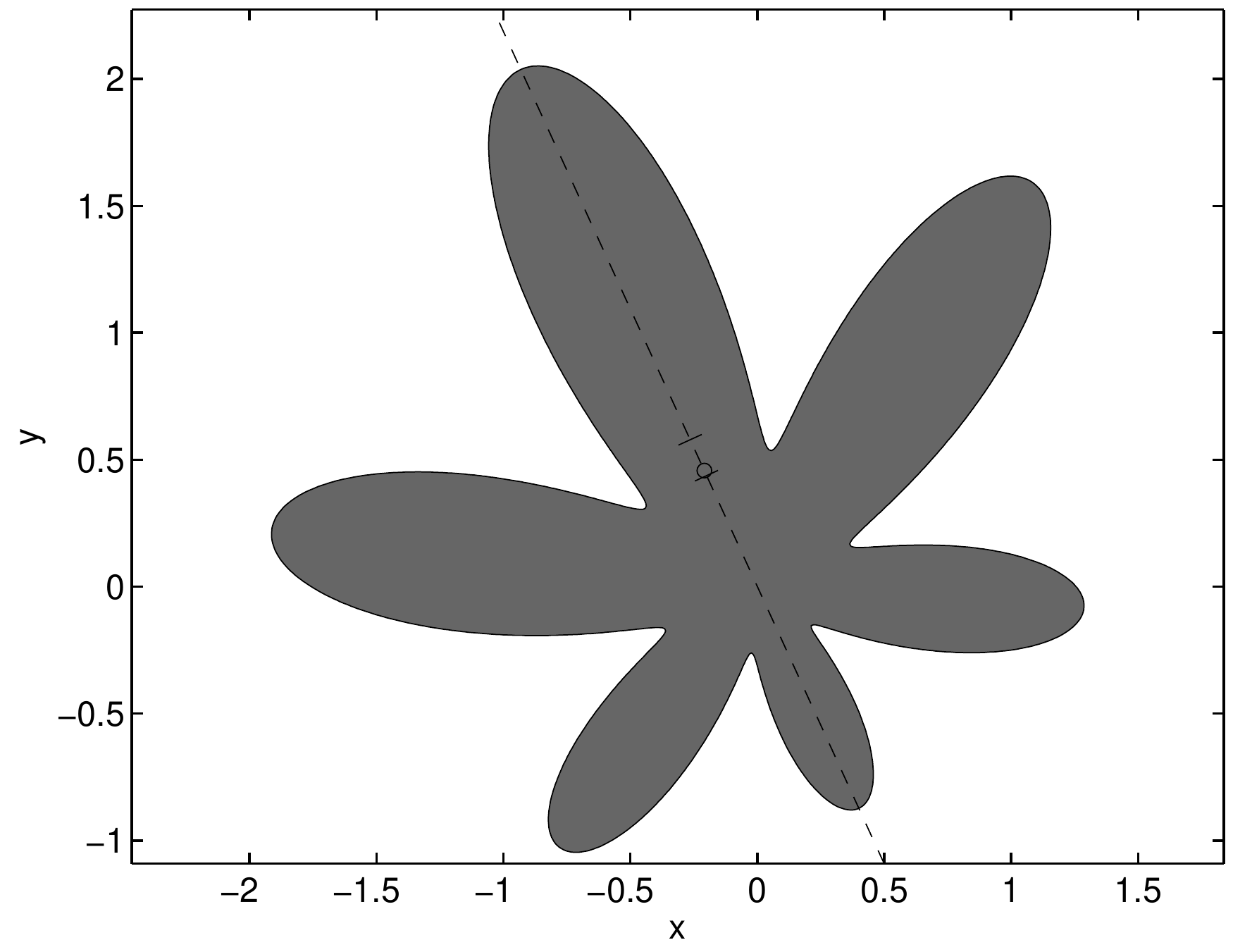}}
\caption{The flower domain. The dashed line is the axis of symmetry and the circle represents the center of gravity of the drop. The interval depicted along the axis of symmetry is where the steady state circle centers end up.}
\label{fig:slanted}
\end{figure}
\begin{figure}[h]
\centering
\resizebox{60mm}{!}{\includegraphics{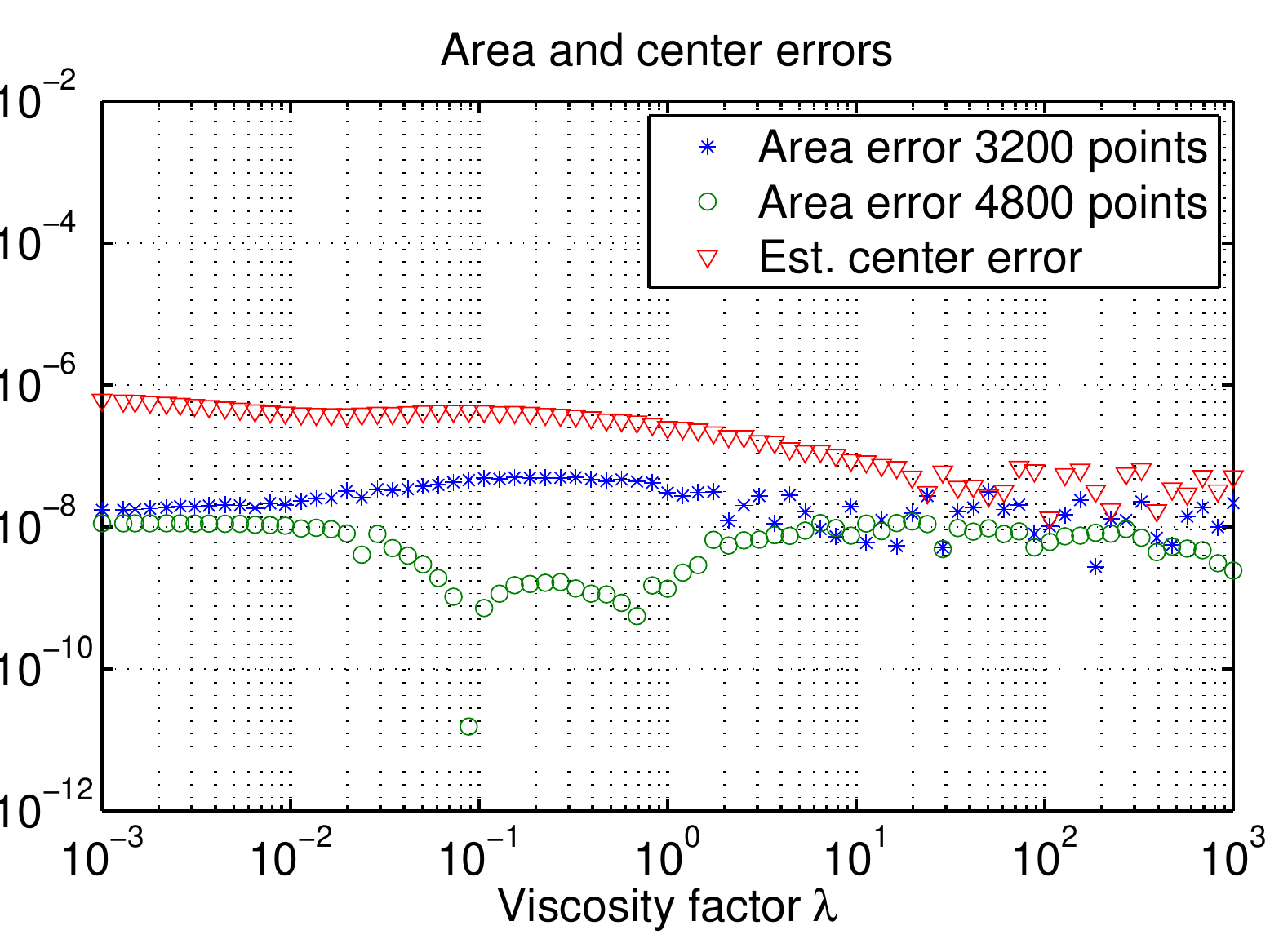}}
\resizebox{60mm}{!}{\includegraphics{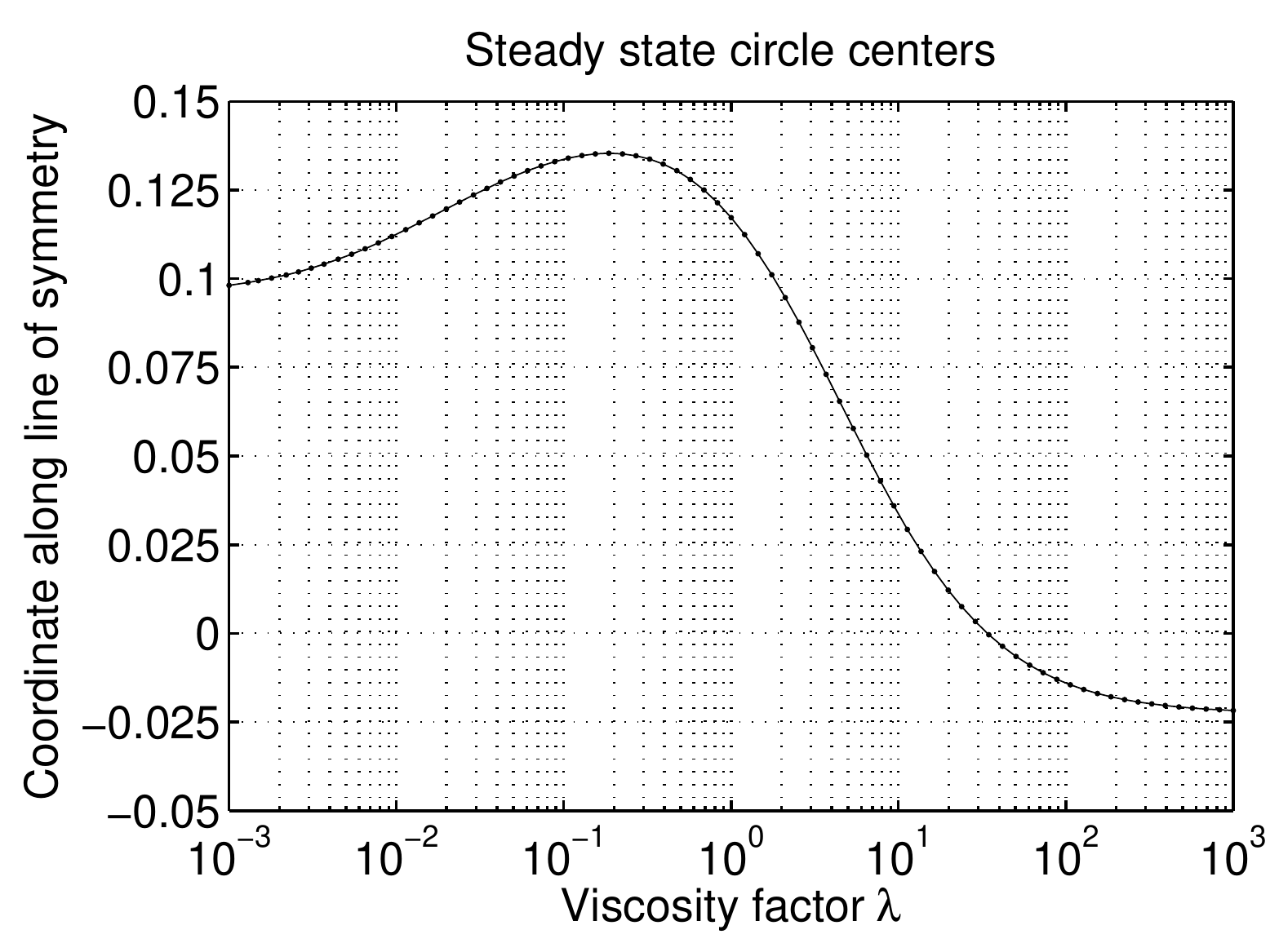}}
\caption{Estimated errors and steady state circle positions for the asymmetric flower. Left : the area and steady state circle center errors as a function of the viscosity ratio $\lambda$. Right : the coordinate, $z_{\rm center}$, of the steady state circle center along the line of symmetry as a function of $\lambda$.}
\label{fig:slerrs}
\end{figure}

The flower drop, depicted in figure~\ref{fig:slanted}, is parameterized as
\begin{equation}
z(s) = e^{{\rm i}(s+2)}\left(1+0.6\cos(6s)\right)\left(1+0.4\cos(s)\right), \qquad s \in [0,2\pi].
\end{equation}
Its shape will cause the center of gravity of the drop to move along the axis of symmetry during the simulation. The total distance it moves depends in a non-trivial manner on the viscosity ratio $\lambda$, as can be seen in figure \ref{fig:slerrs} and table \ref{tab:slanted}. 
Admittedly, for this drop shape the need for special quadrature treatment is not great, but the focus here is instead to study a wide range of drop viscosity ratios and how the centers of the steady state circles vary with $\lambda$.

\begin{table}
\centering  

\begin{tabular}{c c c c c} % centered columns (4 columns)
\hline\hline                        %inserts double horizontal lines
$\lambda$ & $A_{\rm err}$ & $c_{\rm err}$ & $z_{\rm center}$ & $t_{\rm steady}$  \\ [0.5ex] % inserts table 
%heading
\hline                  % inserts single horizontal line
0.001 & 1.7e-8 & 6.1e-7 & (-0.250042,0.546352) & 4.77  \\ 
0.01 & 2.3e-8 & 4.0e-7 & (-0.256018,0.559410) & 4.90  \\
0.1 & 3.0e-8 & 4.3e-7 & (-0.264824,0.578650) & 5.79 \\ 
1 & 3.0e-8 & 2.5e-7 & (-0.257990,0.563718) & 11.3 \\ 
10 & 1.4e-8 & 8.6e-8 & (-0.2232233,0.4877517) & 53.6 \\ 
100 & 9.2e-9 & 1.2e-8 & (-0.2033712,0.4443741) & 458 \\ 
1000 & 2.1e-8 & 5.1e-8 & (-0.2001502,0.4373362) & 4500 \\ 
[1ex]      % [1ex] adds vertical space
\hline %inserts single line
\end{tabular}
\caption{Results from seven simulations on the flower drop with different $\lambda$. $A_{\rm err}$ is the area error, $c_{\rm err}$ is the estimated error in the position of the steady state circle center, $z_{\rm center}$ is the steady state circle center position and $t_{\rm steady}$ is the approximate number of time units needed to reach steady state.}
\label{tab:slanted}
\end{table}

We use an initial number of 3200 discretization points in the simulation for all values of $\lambda$, and 4800 discretization points are used to compute the reference solution. Towards the end of the simulation the number of points has dropped to 1408 and 2096, respectively, due to the adaptivity scheme that reduces the number of discretization points as the circumference of the drop shrinks. The number of discretization points required for sufficient resolution of the drop boundary depends on the velocities and forces due to the surface tension. These quantities are controlled by $\lambda$, and for $\lambda$ in the neighborhood of 1 the number of discretization points required is significantly lower than 3200. For $\lambda = 1$, for example, we only need around 1500 discretization points to reach an area error of $10^{-8}$. Ideally, some a priori refinement taking into account the value of $\lambda$ and the shape of the drop boundary to create an optimal or near-optimal initial grid could be used. Constructing such a scheme seems very difficult, though, and for simplicity we stick with a constant initial grid size.

As can be seen to the left in figure \ref{fig:slerrs}, all the area errors are around $10^{-8}$ in magnitude, which is the set tolerance of the time stepper. This indicates that the spatial resolution is sufficient and that it is the temporal resolution that is the limiting factor. The estimated center errors are somewhat higher, especially for low values of $\lambda$, although never higher than $10^{-6}$. Furthermore, the center errors are stable even though the $\lambda$-values span six orders of magnitude.

To the right in figure \ref{fig:slerrs} we show the coordinate of the circle center $z_{\rm center}$ at steady state along the line of symmetry. We have set $z_{\rm center} = 0$ at the center of gravity of the initial drop and positive coordinates indicates that $z_{\rm center}$ is above and to the left of the center of gravity along the line of symmetry. The value of $z_{\rm center}$ as a function of $\lambda$ is not monotone and there is a maximum for $\lambda \approx 0.2$.

The wall clock time for a simulation run varies greatly, from 1 minute for $\lambda = 1$ to 40 minutes for $\lambda = 10^{-3}$. For high values of $\lambda$ the time taken is around 15 minutes. The reason for the longer simulation times for low values of $\lambda$ is that the number of GMRES iterations required to solve (\ref{eq:disceq}) is higher.

\subsection{The C-domain}
\label{sec:cshape}
\begin{figure}[h]
\centering
\resizebox{85mm}{!}{\includegraphics{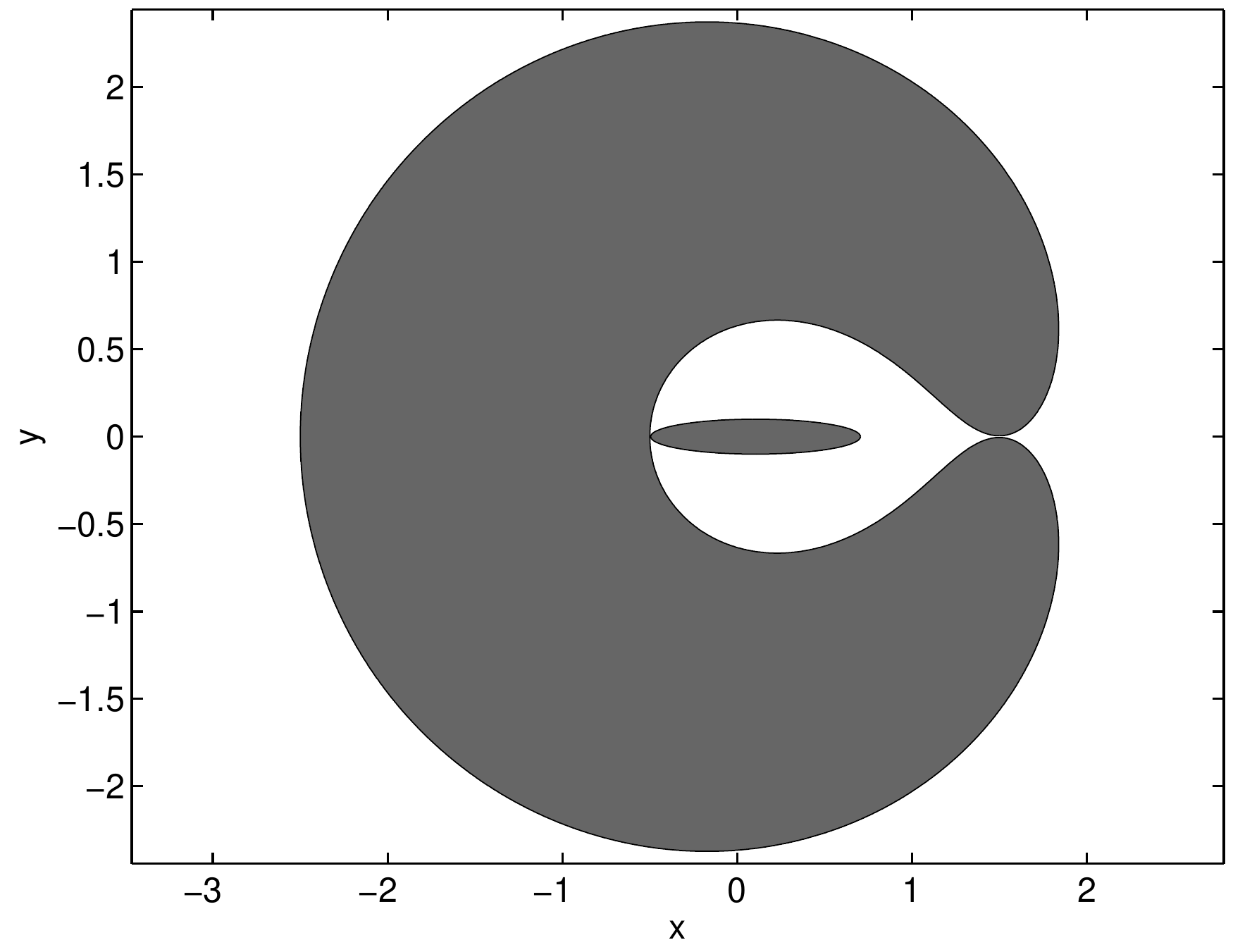}}
\caption{The C-domain.}
\label{fig:cshape}
\end{figure}
The C-domain is a test case with an initial geometry as shown in figure \ref{fig:cshape}
consisting of a C-shaped domain with parameterization
\begin{equation}
z_c(s) = -\left(1.5+\sin\left(s\right)\right)e^{-0.999{\rm i}\pi\cos(s)}, \qquad s \in [0,2\pi],
\end{equation}
and an ellipse described by
\begin{equation}
z_e(s) = 0.6\cos(s)+0.1{\rm i}\sin(s)+0.105, \qquad s \in [0,2\pi].
\end{equation}
The distance between the tips of the ''C'' is 0.001 and the distance between the ellipse and the ''C'' is roughly 0.005, which is close enough to require special quadrature treatment unless the discretization is made very dense. Assigning different values of $\lambda$ to the two drops, we get quite different behavior during the simulation. When the ellipse has a high value of $\lambda$ and the ''C'' a low value, for example, the ellipse tugs at the boundary of the ''C'' quite a lot. This local but aggressive interaction puts some demand on the resolution of the discretization. Some local adaptivity would be very advantageous, but lacking such machinery, we have to make do with increasing the level of refinement on the whole boundary. We use 4800 discretization points on the ''C'' because of this and 800 on the ellipse. At the end of the simulation run these numbers have dropped to 3312 and 512, respectively, due to the adaptivity scheme. The reference simulation is again done with 50\% more points. 

\begin{figure}[htbp]
\centering
\resizebox{90mm}{!}{\includegraphics{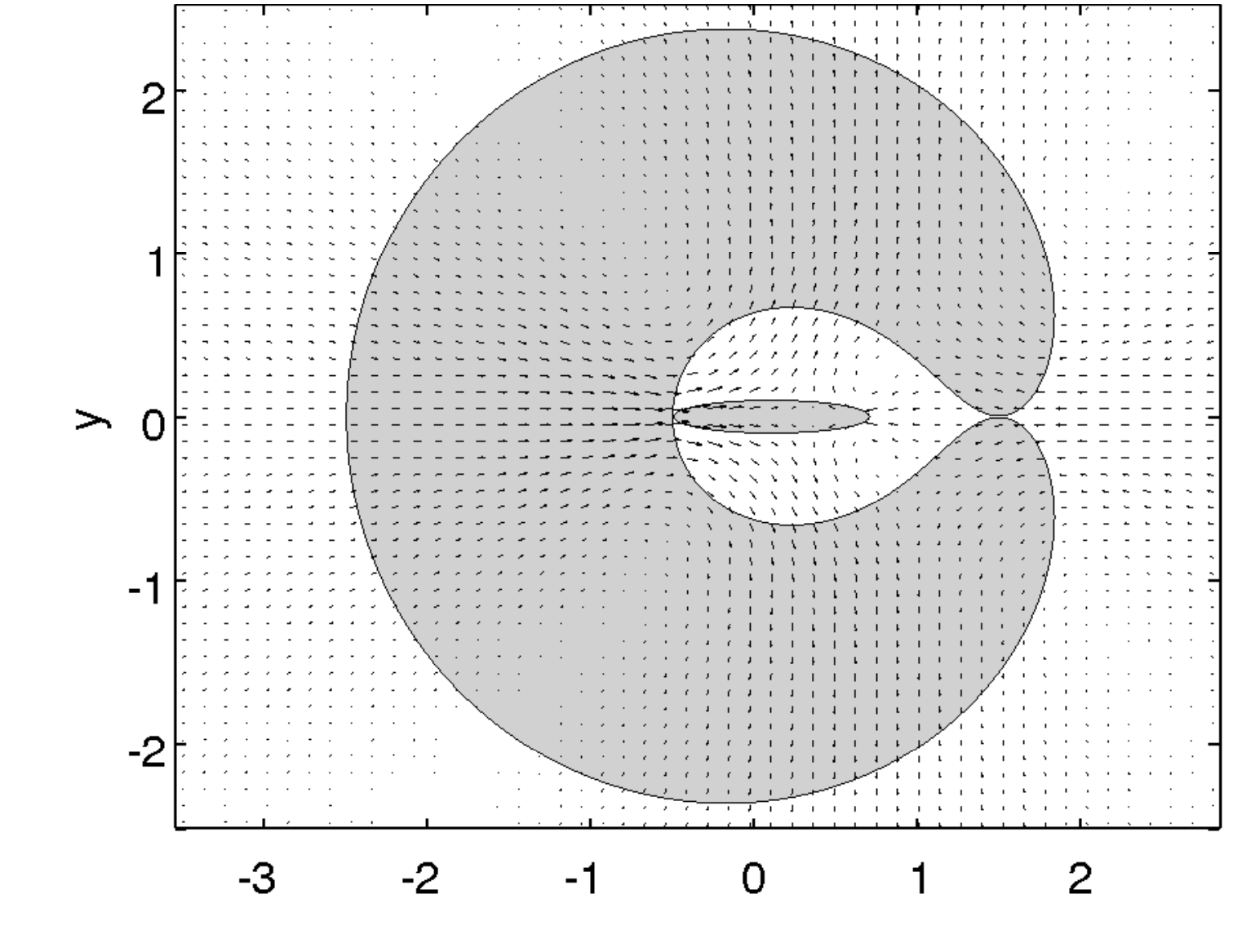}}
\resizebox{90mm}{!}{\includegraphics{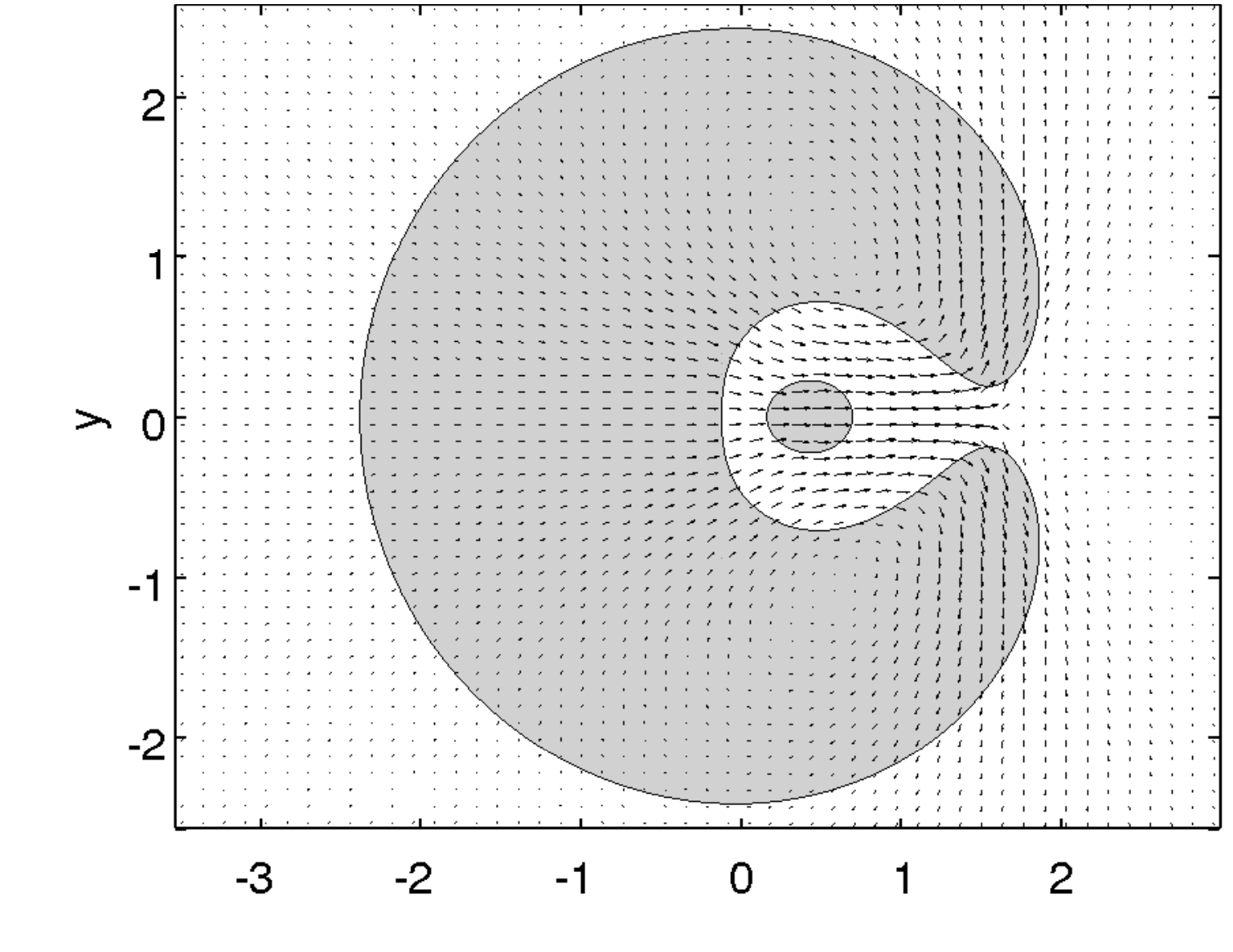}}

\caption{The velocity field of the C-domain with $\lambda^{(1)}=\lambda^{(2)}=1$, at simulation times $t=0$ and $t=3$. The special quadrature is used for field points close to a drop boundary.}
\label{fig:cshapevel}
\end{figure}

In figure \ref{fig:cshapevel} we show quiver plots of the velocity field at $t=0$ and $t=3$ for $\lambda^{(1)}=\lambda^{(2)}=1$. The velocity is computed using (\ref{eq:vel}), and for $z$ close to a drop boundary we utilize special quadrature for accuracy.

\begin{table}
\centering
\begin{tabular}{c c c c c c c c} % centered columns (4 columns)
\hline\hline                        %inserts double horizontal lines
$\lambda^{(1)}$ & $\lambda^{(2)}$ & $A_{\rm err}$ & $c_{\rm err}^{(1)}$ & $c_{\rm err}^{(2)}$ &  $x_{\rm center}^{(1)}$ & $x_{\rm center}^{(2)}$ &$t_{\rm steady}$  \\ [0.5ex] % inserts table 
%heading
\hline                  % inserts single horizontal line
0.1 & 0.1 & 4.6e-9 & 1.1e-5 & 4.1e-7 & 0.021777 & 2.85817 & 15.6  \\ 
0.1 & 1 & 5.6e-10 & 7.3e-6 & 7.0e-6 & 0.024304 & 2.91348 & 15.7  \\
0.1 & 10 & 3.4e-9 & 7.7e-6 & 2.8e-7 & 0.038306 & 2.78027 & 19.9  \\
1 & 0.1 & 1.6e-9 & 1.5e-7 & 3.7e-7 & -0.1123831 & 2.666463 & 31.2  \\
1 & 1 & 1.1e-9 & 9.6e-8 & 3.5e-7 & -0.1107529 & 2.724521 & 31.2  \\
1 & 10 & 1.6e-9 & 5.1e-6 & 7.9e-6 & -0.116858 & 2.72017 & 31.1  \\
10 & 0.1 & 1.4e-9 & 1.7e-7 & 1.2e-6 & -0.321246 & 2.34729 & 155  \\
10 & 1 & 5.5e-9 & 3.1e-7 & 1.1e-6 & -0.320471 & 2.40128 & 154  \\
10 & 10 & 2.5e-9 & 1.9e-8 & 4.9e-8 & -0.3253772 & 2.443563 & 154  \\
[1ex]      % [1ex] adds vertical space
\hline %inserts single line
\end{tabular}
\caption{Results from simulations on the C-shape drop configurations. The superscripts 1 and 2 denotes the large and small drop respectively. $A_{\rm err}$ is the area error, $c_{\rm err}$ is the estimated error in the position of the steady state circle center, $x_{\rm center}$ is the $x$-coordinate of the steady state circle center position and $t_{\rm steady}$ is the approximate number of time units needed to reach steady state.}
\label{tab:cshape}
\end{table}

\begin{figure}[htbp]
\centering
\resizebox{70mm}{!}{\includegraphics{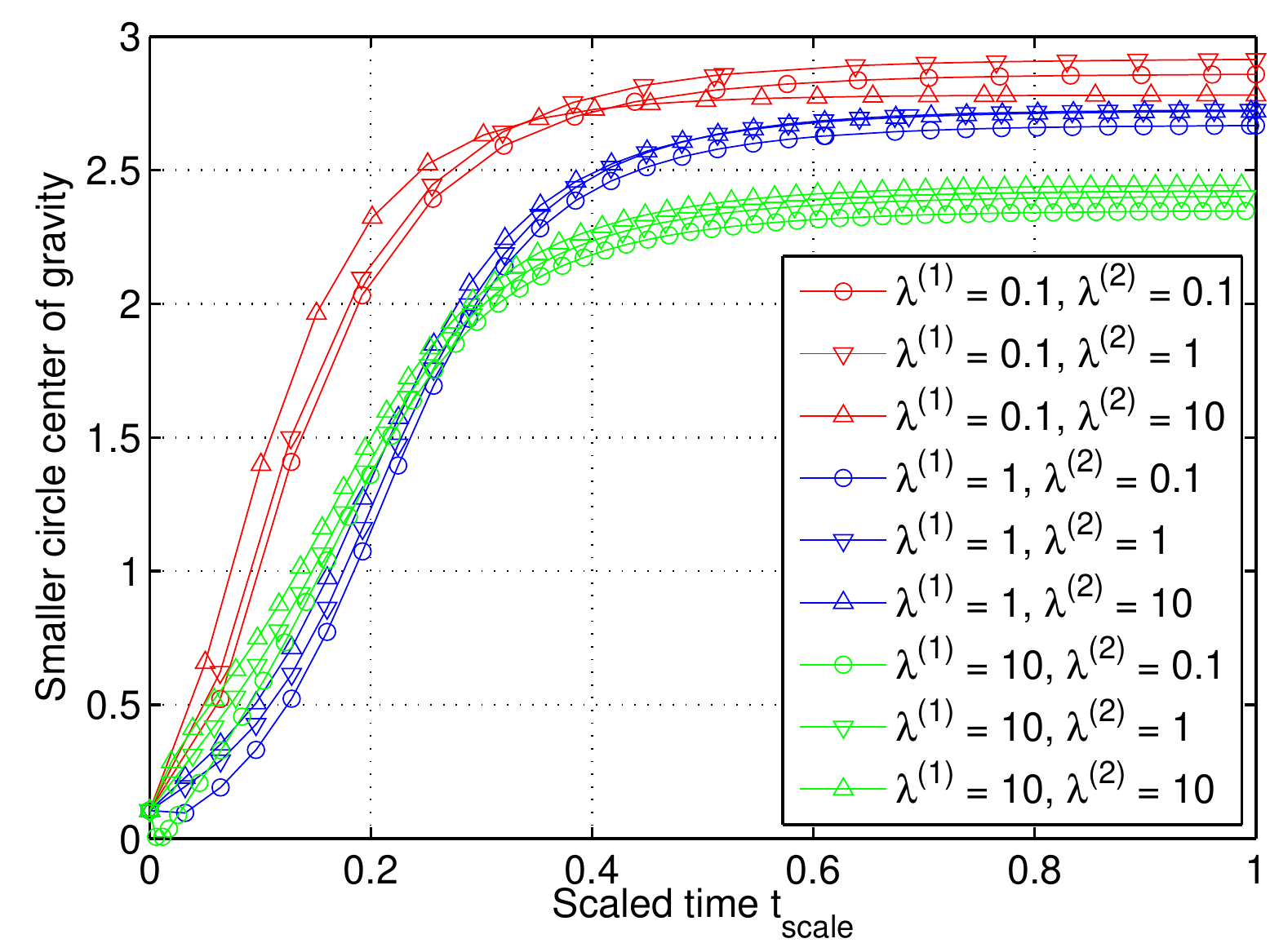}}
\caption{The $x$-coordinate of the center of gravity of the smaller drop in the C-shape drop configuration as a function of normalized time. As is displayed in table \ref{tab:cshape}, the time to steady state depends strongly on the viscosity of the larger drop. Therefore, the time variable has been normalized to $t_{\rm scale} \in [0,1]$ for all runs to clarify the behavior.}
\label{fig:cshapex}
\end{figure}

In table \ref{tab:cshape} we have gathered some simulation results for different drop viscosities. We only show the $x$-coordinate of the steady state circle center position. The final center positions in the $y$ direction are expected to be $0$ due to the symmetry of the drop configuration, and the deviations from $0$ are in all cases less than the center errors in table \ref{tab:cshape}. Furthermore, we see in table \ref{tab:cshape} that the area is preserved to high accuracy, and again that the errors in the steady state circle positions tend to be somewhat higher. There is some indication that for higher values of $\lambda$, with in some sense stiffer drops, the errors are lower here as well. The behavior of the drops seems to be in some sense controlled by the higher viscosity drop: the lower viscosity drop is to a greater degree affected by the velocity field induced by the surface tension of the other drop. The approximate time to steady state depends strongly on the viscosity of the larger drop, as its larger mass takes longer to shift. 

In figure \ref{fig:cshapex}, we show the $x$-coordinate of the center of gravity of the smaller drop over time during the simulation for different combinations of viscosities. The behavior of the smaller drop varies quite a lot in the different cases, but it is especially sensitive to the viscosity of the larger drop. This is certainly so for low viscosities for the smaller drop : for $\lambda^{(1)} = 10, \lambda^{(2)} = 0.1$, and to some degree for $\lambda^{(1)} = 1, \lambda^{(2)} = 0.1$, the center of gravity of the smaller drop actually travels left in the initial stages of the simulation. 

\subsection{The Swiss roll}
\begin{figure}[htbp]
\centering
\resizebox{85mm}{!}{\includegraphics{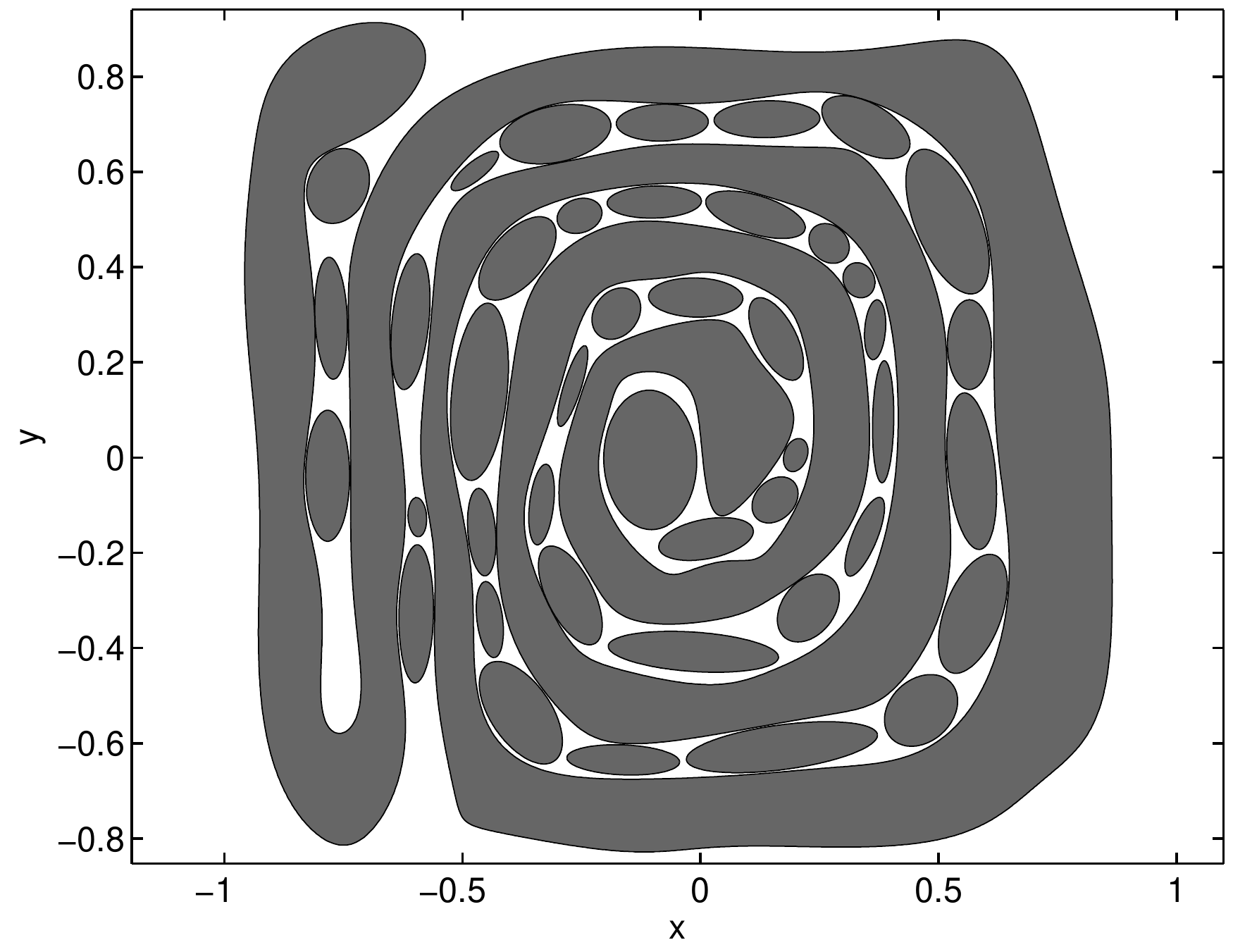}}
\caption{The Swiss-roll domain.}
\label{fig:swissroll}
\end{figure}
We conclude the numerical experiments with a simulation on a rather complicated drop set up. It consists of a spiral and 43 ellipses, and most of the latter are very close to the spiral. At first glance it is perhaps hard to envision the spiral being able to retract itself into a circle with all the ellipses being in the way, but this indeed what happens. The flow field induced by the surface tension is very complicated, especially when the viscosities differ by a great deal, so we will only show results for the case when all $\lambda_k = 1$. Simulations with $\lambda_k = 10$ for the ellipses have been tested, but the spatial resolution needed to resolve the flow field combined with the fact that each velocity calculation requires the solution of a rather large system of equations makes this a daunting task for a mere off the shelf workstation. A full simulation to steady-state would take days of computation.

\begin{figure}[htbp]
\centering
\resizebox{60mm}{!}{\includegraphics{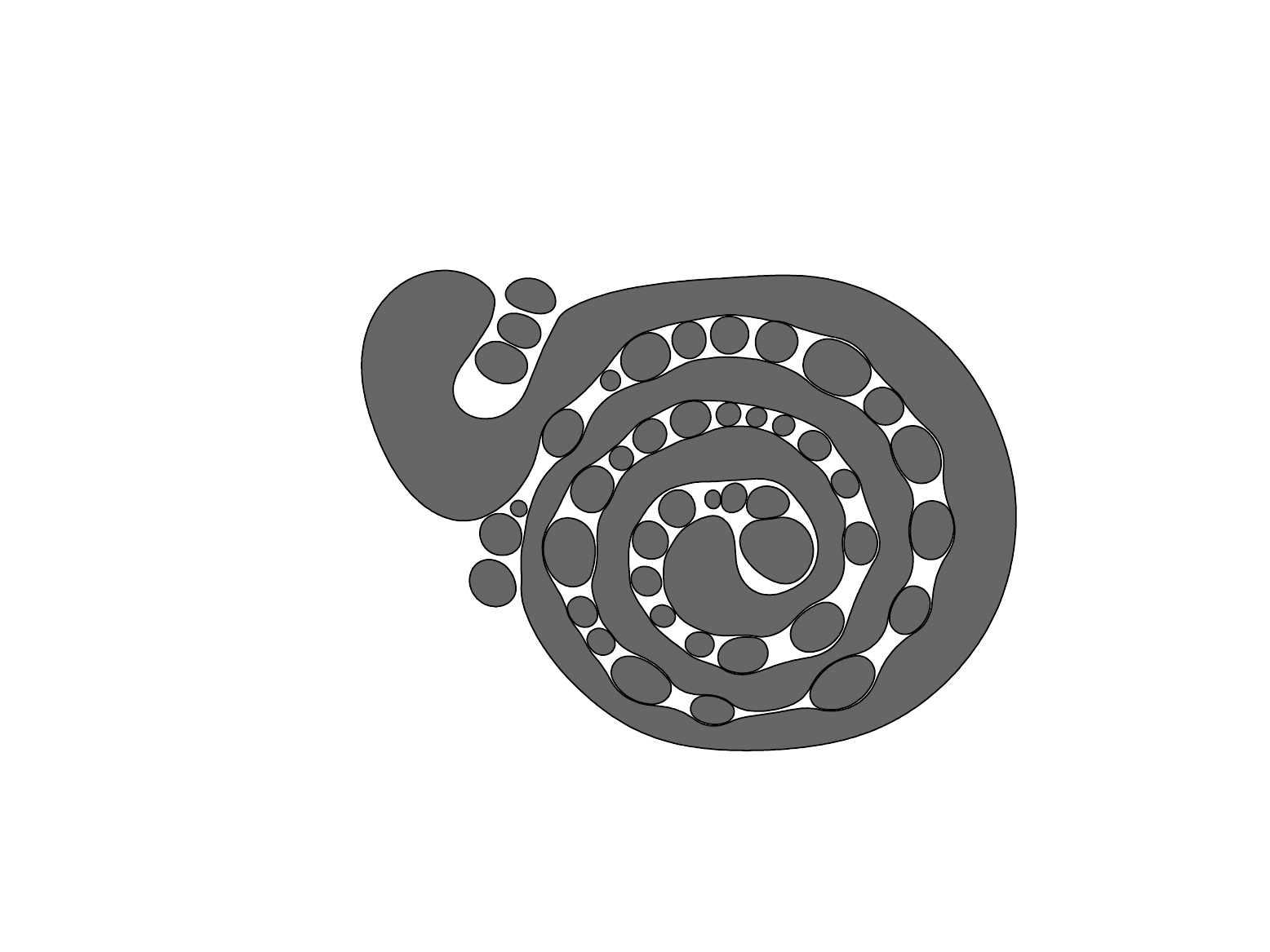}}
\resizebox{60mm}{!}{\includegraphics{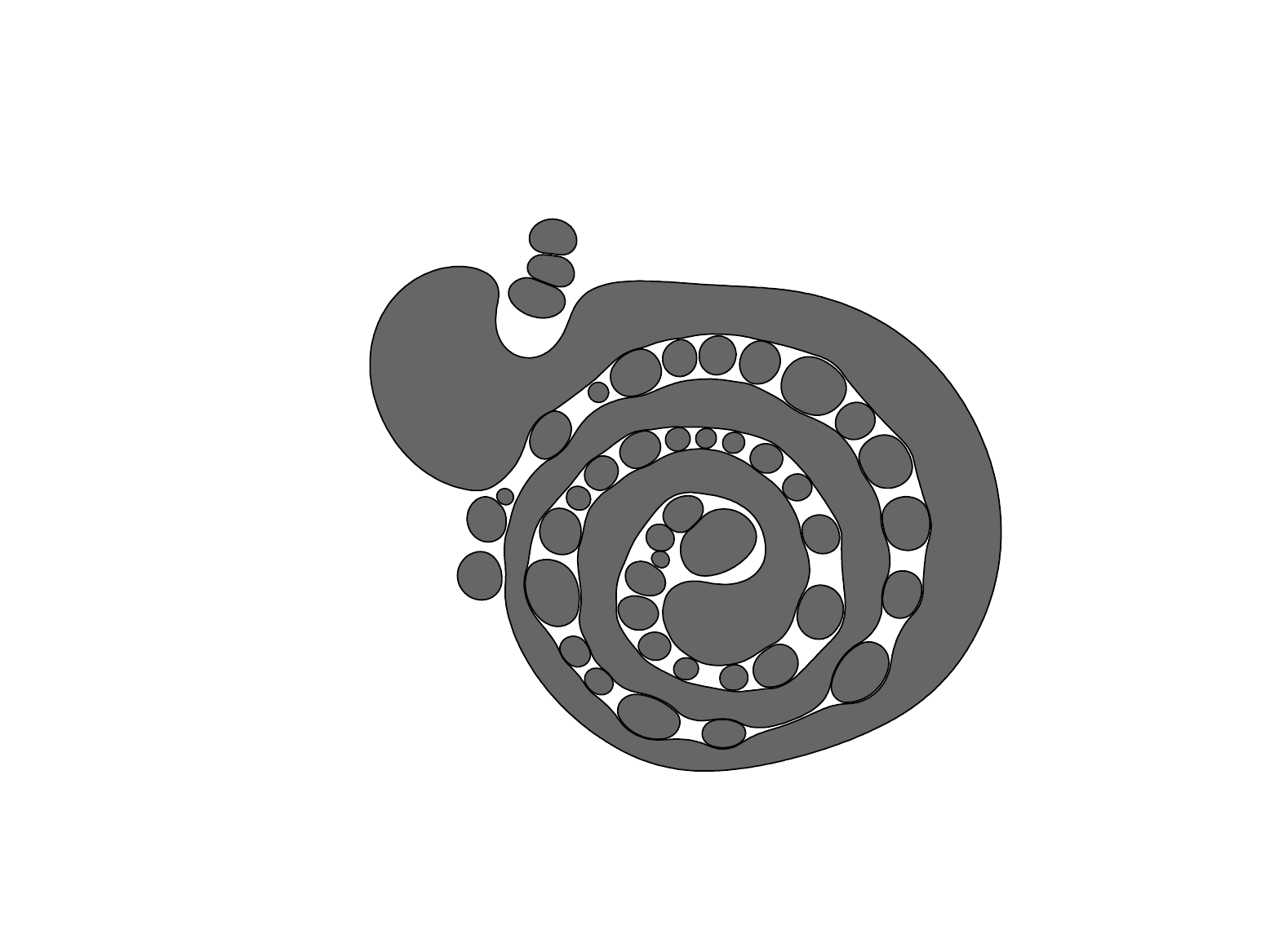}}
\resizebox{60mm}{!}{\includegraphics{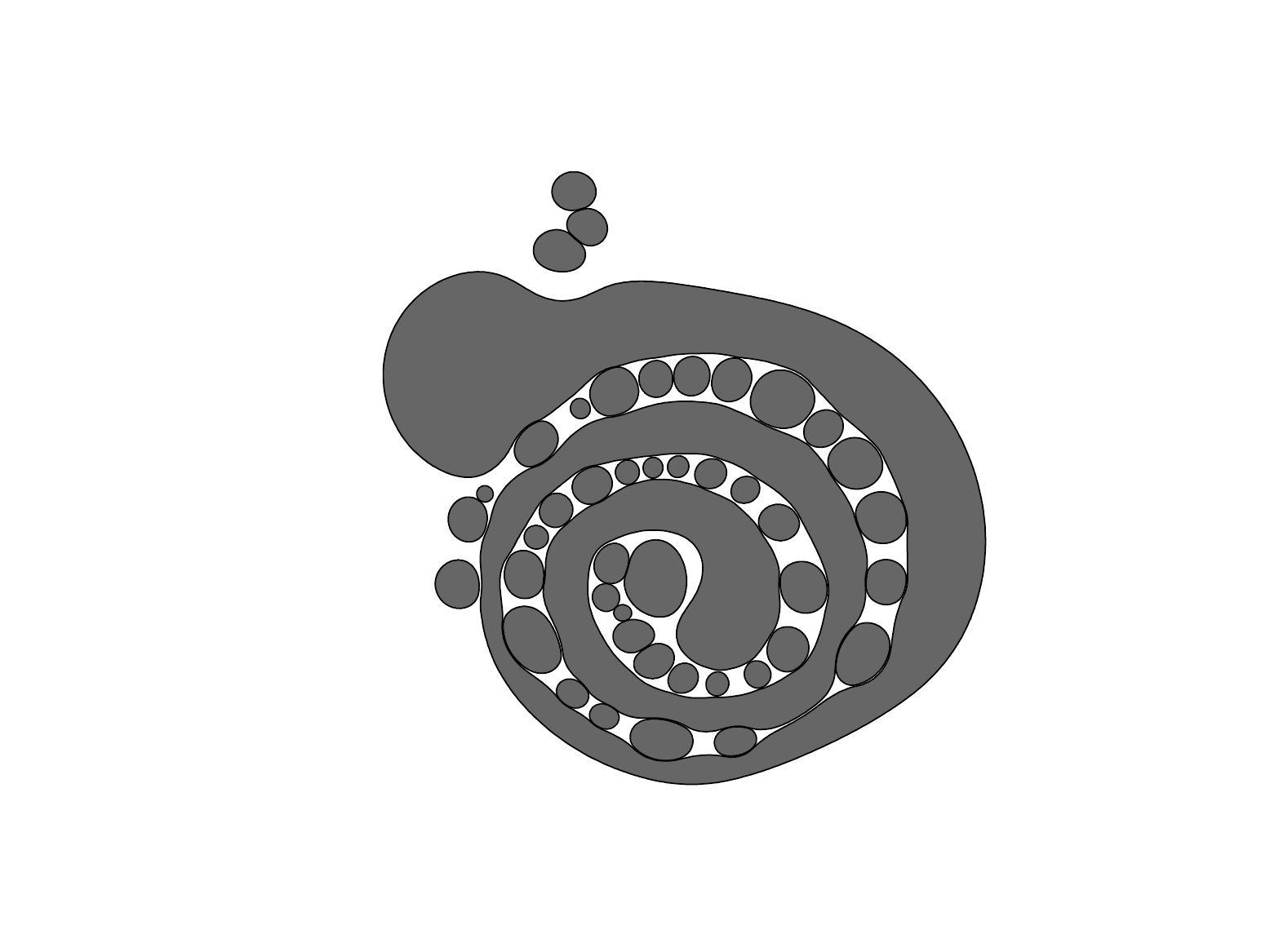}}
\resizebox{60mm}{!}{\includegraphics{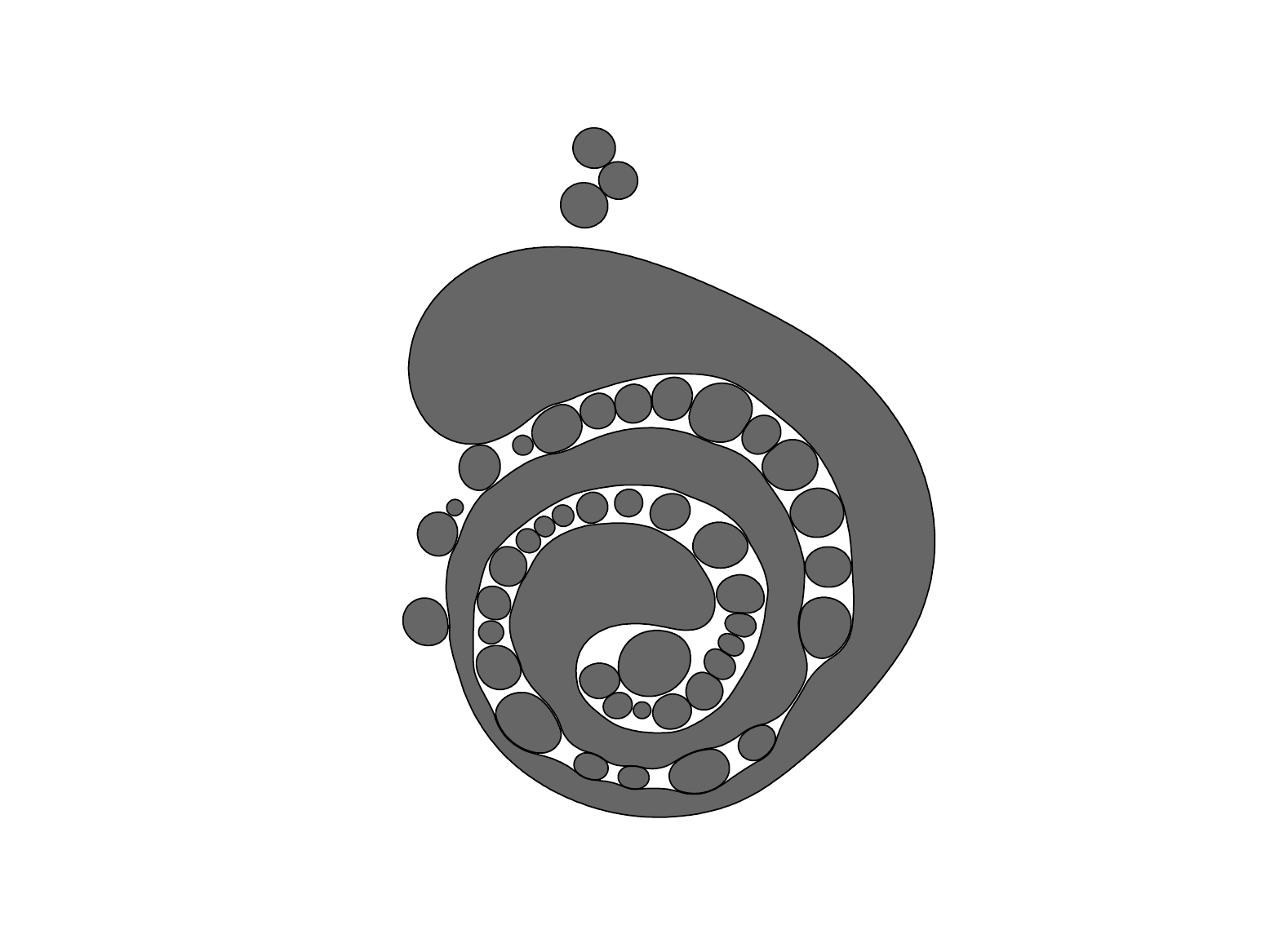}}
\resizebox{60mm}{!}{\includegraphics{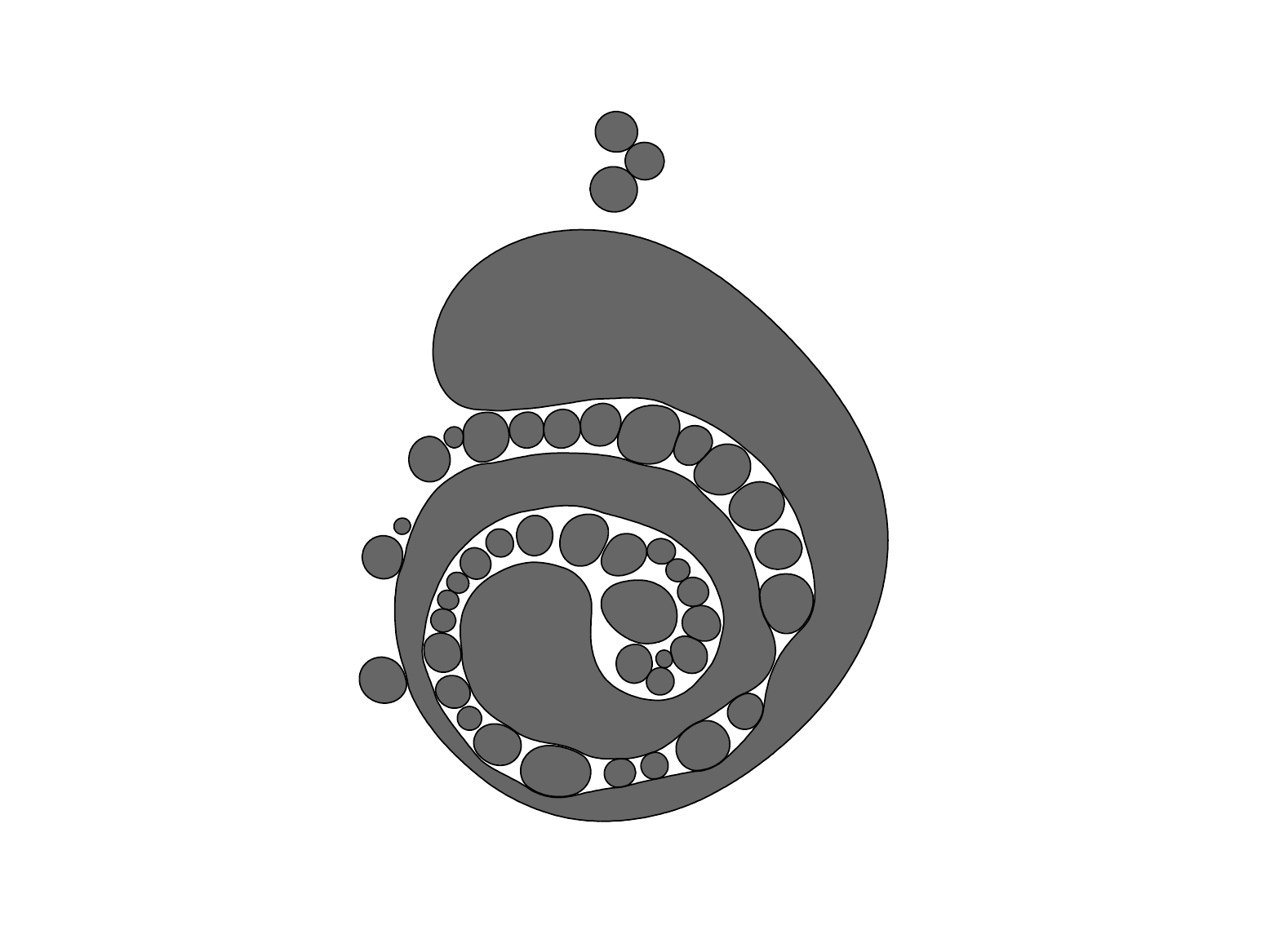}}
\resizebox{60mm}{!}{\includegraphics{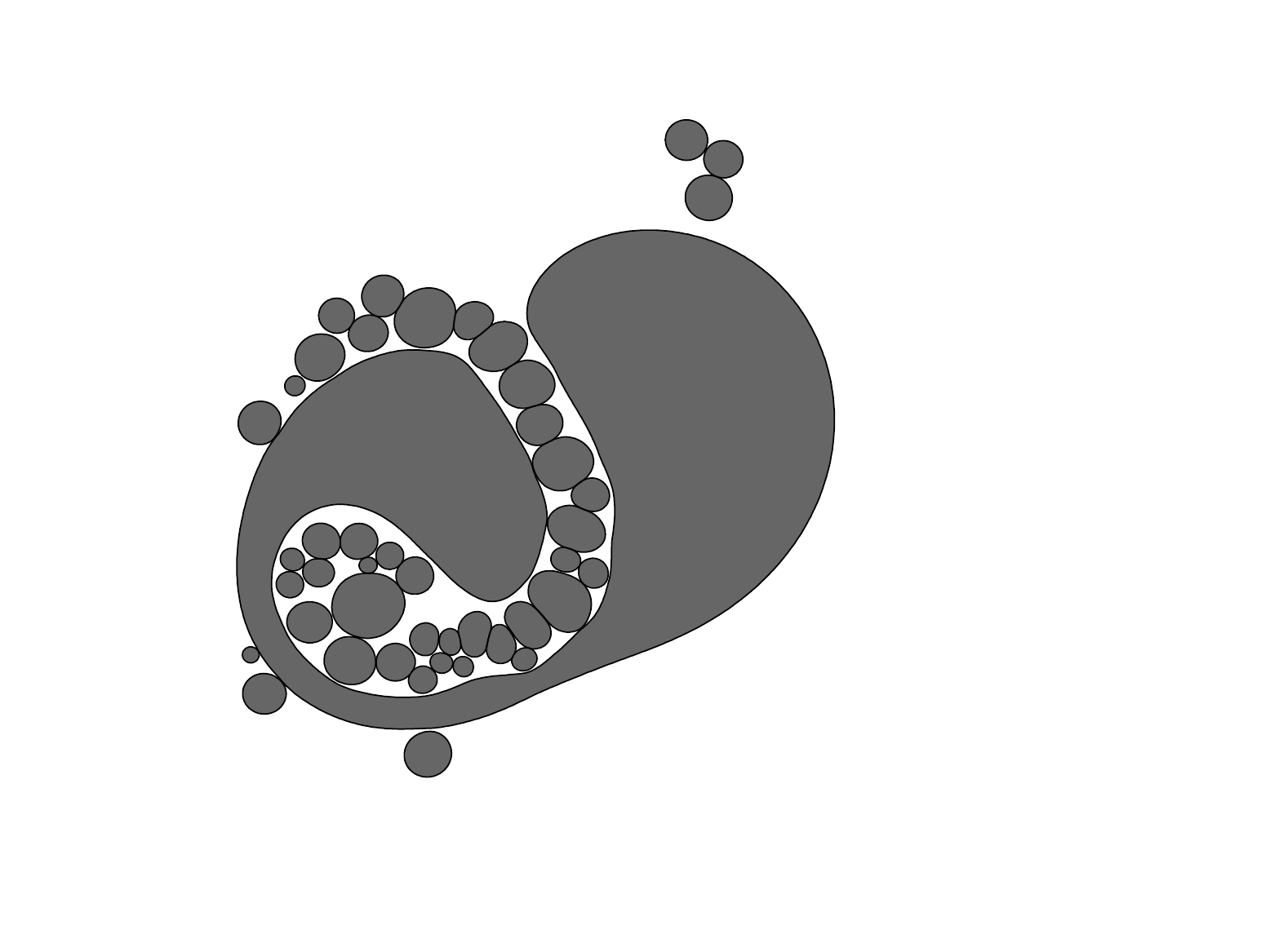}}
\resizebox{60mm}{!}{\includegraphics{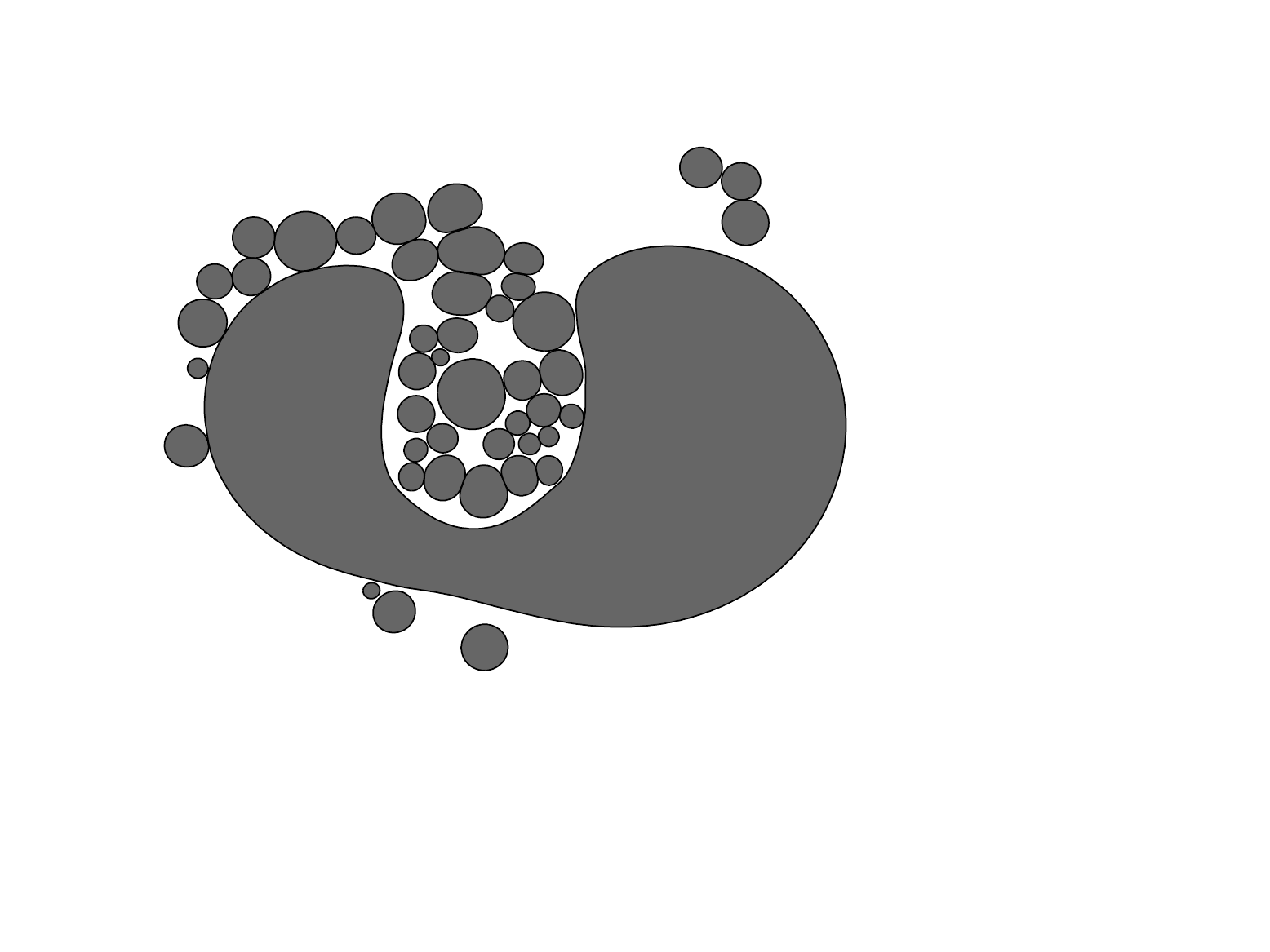}}
\resizebox{60mm}{!}{\includegraphics{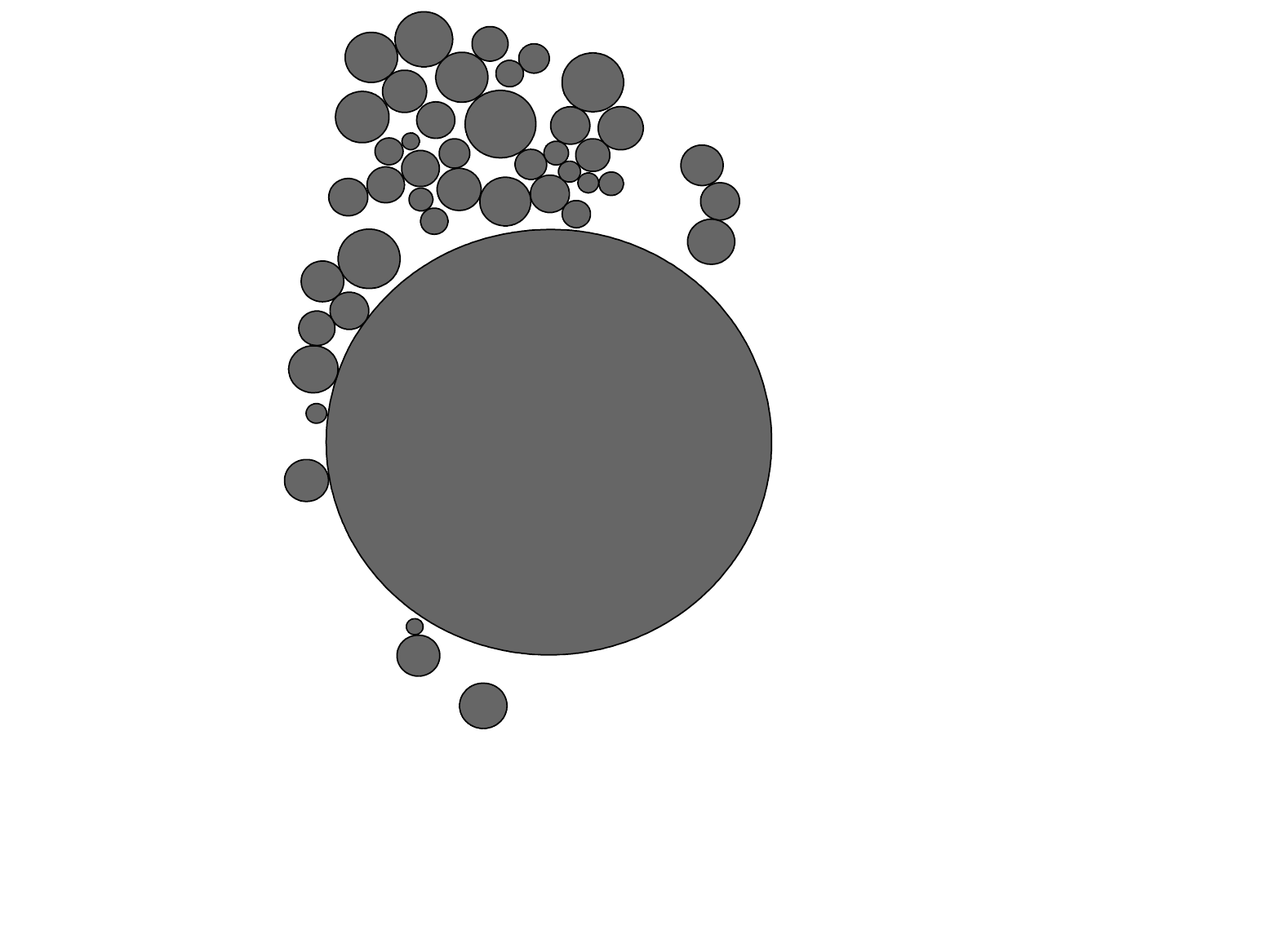}}
\caption{The evolution of the Swiss-roll drop configuration. From the top-left to the bottom right the simulation times are : $t=2.25, 3.4, 4.5, 6.75, 9, 18, 27, 45$.}
\label{fig:swissrollsim}
\end{figure}

\begin{figure}[htbp]
\centering
\resizebox{90mm}{!}{\includegraphics{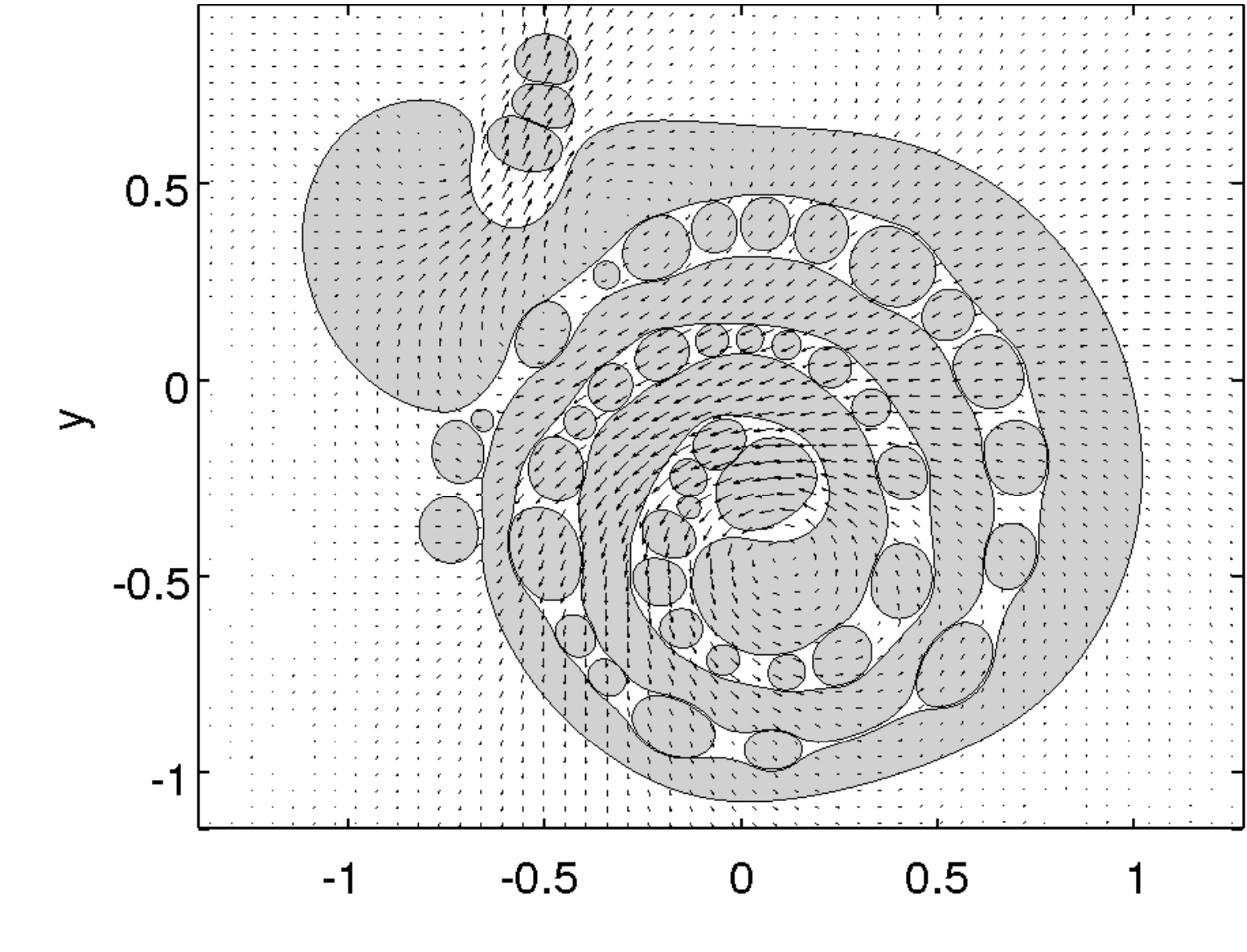}}
\resizebox{90mm}{!}{\includegraphics{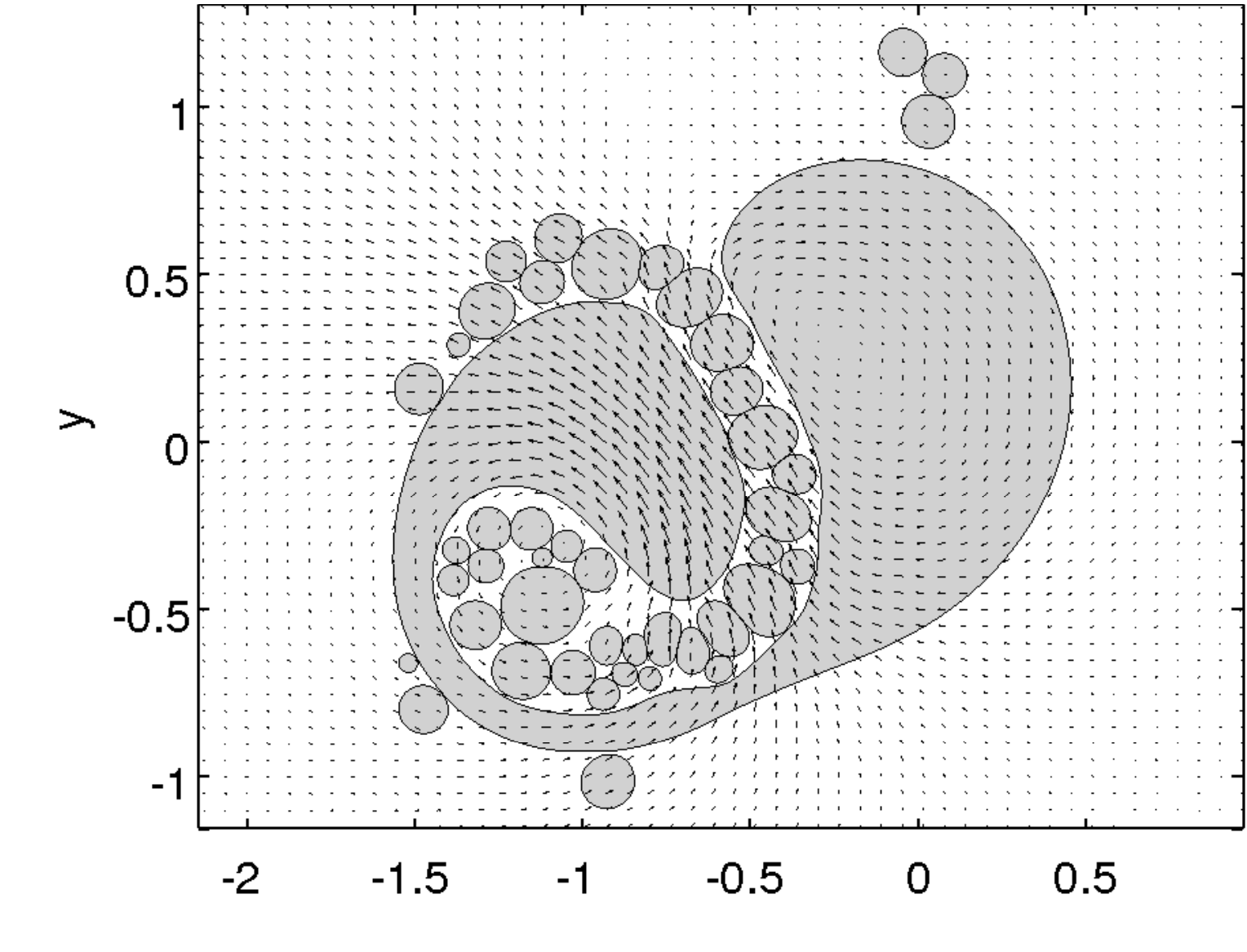}}

\caption{The velocity field of the Swiss roll at simulation times $t=3.4$ and $t=18$. The special quadrature is used for field points close to a drop boundary.}
\label{fig:swissvel}
\end{figure}
We use 25200 discretization points for the boundaries, 8000 for the spiral and 400 for each ellipse. As before, we compute a reference solution with 50\% more points. This is a low density compared to the examples above, but since we are only concerned with $\lambda_k = 1$ here, it suffices. The entire computation to approximate steady state takes about 8 hours. The area error is on the order of the Runge-Kutta tolerance $10^{-8}$ as before, and the maximum circle center error is $1.4\cdot 10^{-5}$ which is a bit higher than for the simpler configurations in the previous examples. Handling a domain of this complexity with a reasonable degree of accuracy would be infeasible without using special quadrature : the number of discretization points needed to resolve the lubrication effects would be enormous. Figure \ref{fig:swissrollsim} shows a number of frames from the evolution of the swiss roll drop configuration, while figure \ref{fig:swissvel} shows the velocity field at two instances during the simulation. 
%The full 25 fps animation in animated .gif format can be downloaded from %\texttt{www.nada.kth.se/\textasciitilde annak/}.

\section{Conclusions and outlook}
\label{sec:outlook}
We have presented methods to accurately evolve the boundaries of general setups of bubbles in quasi-static two dimensional Stokes flow. The main novelty of the paper is the introduction of a  specialized quadrature scheme enabling the treatment of significantly more complicated geometries than was previously possible. Furthermore, the high local accuracy achievable with the specialized quadrature scheme allows for correct modeling of for example lubrication and mass preservation without the need for artificial constraints. This also keeps the solver cleaner and simpler. The capabilities of the solver are demonstrated via a number of numerical experiments and several benchmark results are reported to facilitate validation and comparisons to other solvers. 

There are numerous ways to extend the scope of the solver. For example periodicity in one or more directions and introduction of fixed walls of arbitrary shapes have been considered. This requires a new integral equation formulation, but still allows for using the machinery developed in this paper. The step thereafter is to add additional physics to the problem, such as that of surfactants and electric fields. This is highly relevant to droplet based micro-fluidic system where surfactants are used for stabilization and electric fields for sorting. Using primitive variable formulations it is also possible to extend to surfaces in 3D, but it remains to be seen what efficiency can be achieved. There are indications \cite{Helsing2013} that the special quadrature approach used here can be generalized to surfaces in 3D, and there are also other promising approaches such as the aforementioned QBX method by Kl\"ockner, et al. \cite{Kloeckner2013}.

%%\newpage

\bibliography{stokes_short}
\bibliographystyle{plain}

\end{document}